\documentclass[final]{siamltex}

\usepackage{cleveref}
\usepackage{algorithm}
\usepackage{enumerate}
\usepackage[left]{showlabels}
\usepackage{shorttoc}
\usepackage{siunitx}
\usepackage{bm}
\usepackage{amssymb,amstext,mathrsfs}
\usepackage{amsmath}
\usepackage{cite}
\usepackage{float}
\usepackage{graphicx}
\usepackage{multirow}
\usepackage[colorlinks,linkcolor=blue,anchorcolor=black,citecolor=blue]{hyperref}
\usepackage{algpseudocode}
\usepackage{setspace}
\usepackage{simplemargins}

\DeclareMathOperator*{\argmin}{arg\;min}

\numberwithin{equation}{section}
\newtheorem{remark}{Remark}
\newtheorem{example}{{Example}}

\newcommand{\R}{\mathbb{R}}

\title{A Domain-Decomposition Model Reduction Method for Linear convection-diffusion Equations with Random Coefficients\thanks{This material is based upon work supported in part by 
the U.S. Department of Energy, Office of Science, Office of Advanced Scientific Computing Research, Applied Mathematics program under contracts ERKJ259, ERKJ320; the U.S. National Science Foundation, Computational Mathematics program under award 1620027; and by the Laboratory Directed Research and Development program at the Oak Ridge National Laboratory, which is operated by UT-Battelle, LLC., for the U.S.~Department of Energy under Contract DE-AC05-00OR22725.}
}

\author{Lin Mu\thanks{Computer Science and Mathematics Division, Oak Ridge National Laboratory, Oak Ridge, TN 37831 (mul1@ornl.gov, zhangg@ornl.gov).}
\and Guannan Zhang\footnotemark[2]
}

\begin{document}
\maketitle


\begin{abstract}
We develop a domain-decomposition model reduction method for linear steady-state convection-diffusion equations with random coefficients. Of particular interest to this effort are the diffusion equations with random diffusivities, and the convection-dominated transport equations with random velocities. We investigate the equations with two types of random fields, i.e., colored noises and discrete white noises, both of which can lead to high-dimensional parametric dependence. The motivation is to use domain decomposition to exploit low-dimensional structures of local problems in the sub-domains, such that the total number of expensive PDE solves can be greatly reduced. Our objective is to develop an efficient model reduction method to simultaneously handle high-dimensionality and irregular behaviors of the stochastic PDEs under consideration. The advantages of our method lie in three aspects: (i) online-offline decomposition, i.e., the online cost is independent of the size of the triangle mesh; (ii) operator approximation for handling non-affine and high-dimensional random fields; (iii) effective strategy to capture irregular behaviors, e.g., sharp transitions of the PDE solution. Two numerical examples will be provided to demonstrate the advantageous performance of our method.
\end{abstract}

\begin{keywords}
domain decomposition, uncertainty quantification, parametric PDEs, random fields, high dimensionality, sharp transitions
\end{keywords}

\begin{AMS}
65D15, 65N35, 65N12, 65N15, 65C20, 65C30
\end{AMS}

\section{Introduction}
This paper focuses on linear steady-state convection-diffusion equations with random coefficients. Of particular interest to this effort are two types of partial differential equations (PDEs), i.e., the diffusion equations with random diffusivities, and the convection-dominated transport equations with random velocities. We investigate two types of random fields: the colored noises and the discrete white noises. 
The PDEs of interest are widely used to describe subsurface flows in porous media. There are two major challenges in solving such PDEs, i.e.~high dimensional parameterization and irregular behaviors, e.g., solutions with sharp transitions. The parametric dimension depends on the discretization of the random fields. For a colored noise, the most common discretization is the truncated Karhunen-Lo\`{e}ve (KL) expansion \cite{Kac:1947bg}, and the parametric dimension is the number of retained singular values in the truncated expansion. For a discrete white noises, the random field is defined by a piecewise constant function, where the function values for different pieces are independent and identically distributed (i.i.d.) random variables, so that the parametric dimension is the number of pieces in the physical domain. We remark that, in this effort, we only consider {\em finite} dimensional discrete white noises, which is similar to the ``inclusion problem'' considered in \cite{BECK:2014ey} .

In the literature, three common approaches for solving the PDEs of interest are Monte Carlo methods \cite{Fishman:1996es,Peherstorfer:2018dg}, stochastic spectral methods \cite{Cohen:2011jp,Cohen:2010kz,Gunzburger:2014hi,Xiu:2005ia}, and reduced-basis methods \cite{Benner:2015jp,Quarteroni:2015wi,Rozza:2008gg,DeVore:2013eq,Chen:2013gd,Binev:2011fj}. The Monte Carlo methods, including their multilevel/multi-fidelity variants\cite{PWG17MultiSurvey,Giles:2015dd}, are insensitive to the parametric dimension, but feature slow convergence. One particular class of spectral methods utilizes orthogonal polynomials to build sparse approximations of the map from the parameters to the PDE solution. This type of methods have been very successful in exploiting the sparsity of the parametric dependence (e.g., see \cite{Cohen:2010kz}). However, when the dimension is very large and/or the solution map does not have the desired regularity, we will not have sufficient sparsity to build accurate approximations with affordable computational effort. The reduced-basis methods approximate the solution manifold by constructing reduced bases in the finite element space, via proper orthogonal decomposition (POD) or greedy algorithms. In this case, the best convergence rate is described by the decay of the Kolmogorov $n$-width \cite{Cohen:2015ina} with $n$ the dimension of the reduced subspace. Similar to the spectral methods, high-dimensionality and irregularities may lead to a slow decay of the Kolmogorov $n$-width, which will deteriorates the performance of the classic reduced-basis methods. 
 
The domain decomposition (DD) methods were originally proposed to develop parallel solvers for deterministic PDEs (e.g., see \cite{Toselli:2005hi} for details), and the same idea for parallelization can be extended to the stochastic PDE setting \cite{Sarkar:2009eza,Subber:2012gw}. Recently, more research have been conducted on the DD methods for stochastic PDEs, due to the observation that domain decomposition is an effective approach to reduce the parametric dimension. 
For example, according to the results about the dependence of the KL eigenvalues on the ratio of the physcial domain size and the correlation length \cite{Huang:2001fz,Schwab:2006isa,Pranesh:2016dy}, it is easy to see that the KL eigenvalues decay faster as the ratio becomes smaller. A number of attempts using DD methods have being made to alleviate the curse of dimensionality. In \cite{Chen:2015jr}, Chen, et al.~utilized the DD approach to solve stochastic elliptic PDEs, where the local solution in each sub-domain was approximated in a low-dimensional parametric space. In \cite{Hou:2017fu}, Hou, et al.~combine the DD strategy with multi-scale finite element methods to solve elliptic PDEs with high-contrast random medium. In \cite{Contreras:H77YYtfe}, Contreras, et al.~developed a new approach to capture the correlation structure of local random variables in different sub-domains, without performing expensive global KL expansion; and such strategy was incorporated into the DD framework to solve stochastic elliptic PDEs in \cite{Contreras:aZvF515u}. In 
\cite{Tipireddy:2017tz}, Tipireddy, et al.~employed the local KL expansion to reduce local dimensions, and then combined basis adaptation and Hilbert-space KL expansion to approximate local solutions in the sub-domains. However, to our best knowledge, there is still a lack of efficient numerical capability that can simultaneously handle high-dimensionality and irregular behaviors of the stochastic PDEs under consideration; and developing such a capability is the objective of this paper. 

In this effort, we develop a new domain-decomposition model reduction (DDMR) method, which can exploit the low-dimensional structure of local PDE problems from various perspectives. The resulting algorithm can be divided into an offline procedure and an online procedure. The offline procedure consists of four main stages. The first is to divide the physical domain into a set of non-overlapping sub-domains, generate local random fields and establish the correlation structure among local fields. For the discrete white noise, we only need to align the domain partition with the partition of the noise, and there is no correlation among local fields. For the colored noise, we decompose the covariance function in the global domain and sub-domains, respectively, so as to generate the global and local KL expansions. The correlation structure can be obtained by constructing a map, using the least-squares approach, from the global random variables to the local random variables of each sub-domain. Since reducing the cost of generating KL expansion is not an objective of this effort, our strategy may not be optimally efficient, especially when having a very fine triangle mesh. Nevertheless, this operation can be accelerated by using the parallel algorithm developed in \cite{Contreras:H77YYtfe}. The second stage is to generate two sets of training data, i.e., a set of snapshots of the PDE solutions and another set of snapshots of the local stiffness matrices. The third stage is to use singular value decomposition (SVD) to generate a set of reduced bases for the PDE solution in the sub-domains and on the interfaces, then define reduced local stiffness matrices via multiplying each reduced basis by the corresponding blocks of the local stiffness matrices. The fourth stage is to establish sparse approximations to the entries of the reduced local stiffness matrices in low-dimensional local parametric spaces, which finishes the offline procedure. The online procedure is easy to conduct based on the outputs of the offline procedure. When a new realization of the global random field is generated, we map the global random variables to the local random variables, evaluate the sparse approximations of the reduced local stiffness matrices, assemble the reduced global Schur complement matrix, solve the coefficients of the reduced bases on the interfaces, assemble the reduced local Schur complement matrices and solve the coefficients of the reduced bases in the interior of the sub-domains. 

The advantages and contributions of our method lie in the following three aspects. First, the DDMR method has the online-offline decomposition feature, 
i.e., the online computational cost is {\em independent} of the triangle mesh size. This is achieved by utilizing the generated reduced bases. The 
reduced bases on the interfaces are used to reduce the global Schur complement matrix, and the reduced bases in the sub-domains reduce the size of the 
linear system recovering the local solutions. Moreover, since the sizes of the local stiffness matrices are reduced, the total number of entries that need to be 
approximated by sparse polynomials also becomes {\em independent} of the  mesh size, which is critical to the online-offline decomposition. Second, the 
DDMR method can handle the PDEs of interest with non-affine high-dimensional random coefficients. The challenge caused by non-affine coefficients 
is resolved by approximating the entries of the reduced stiffness matrices. The high-dimensionality is handled by the DD strategy. We would like to point out 
that PDEs with discrete white noises are very difficult problems without using the DD strategy. In fact, the independence and isotropy of the large number of 
local random variables make it hard to exploit sparsity to build polynomial approximations. Nevertheless, the DD strategy makes it much easier to solve, 
even easier than the case of having colored noises. Third, the DDMR method can avoid building sparse approximations to local PDE solutions. This 
property is very important in solving the convection-dominated PDE. For instance, the solution in Example \ref{ex2} has a sharp transition caused by the 
boundary condition, and the transition layer moves as the parameter value changes. This irregular behavior will propagate to the parametric space, so that we have to handle sharp transitions when approximating the solution directly. 
However, the entries of local stiffness matrices are not affected at all by such irregularity, so that we can still achieve the spectral convergence in approximating 
the entries of the reduced stiffness matrices. Moreover, if building sparse approximations to the local solutions (e.g., \cite{Chen:2015jr,Hou:2017fu,Contreras:aZvF515u}), we need to decompose 
a local problem into a set of sub-problems, each of which is equipped with a different boundary condition. As such, the total number of local solutions to be 
approximated is the number of the sub-domains multiplied by the total degrees of freedom on the interfaces. In comparison, our approach only approximates one 
reduced local stiffness matrix for each sub-domain. 

The outline of this paper is as follows. In Section \ref{sec2}, we setup the context of this work by introducing the PDEs of interest and the definitions of the random fields under consideration. In Section \ref{sec3}, we briefly recall the deterministic domain decomposition method, which will be used as the exact model. Our DDMR method will be developed in Section \ref{sec:DD_MR}. In Section \ref{example}, we apply our method to the diffusion equation with random diffusivities and the convection-dominated transport equation with random velocities. Finally, some concluding remarks are given in Section \ref{conc}.


\section{Problem setting}\label{sec2}
%
Let $D\subset \R^d, d=1,2,3$, be a bounded 
domain with Lipschitz continuous
boundary, and $(\Omega, \mathcal{F},\mathbb{ P})$ denote a complete probability space, where $\Omega$ is the sample space, $\mathcal{F}\subseteq 2^{\Omega}$ is a $\sigma$-algebra, and $\mathbb{P}$ is the associated probability measure. 
We consider the following stochastic boundary value problem: find a function $u: \overline{D} \times \Omega\rightarrow \mathbb{R}$, such that it holds $\mathbb{P}$-a.e.~in $\Omega$
\begin{equation}\label{Lauf}
\left\{
\begin{aligned}
-\nabla
 \cdot \big(a(x, \omega) \nabla u(x, \omega) \big) + \bm{b}(x, \omega) \cdot \nabla u (x, \omega)  & = f(x) \;\;\;\mbox{  in  }\; D,\\
u(x, \omega) &=\, w(x) \;\;\; \mbox{  on  }\; \partial D,
\end{aligned}\right.
\end{equation}
where $a(x, \omega)$ and $\bm{b}(x, \omega) := (b_1(x, \omega), \ldots, b_d(x, \omega))$ are the random diffusivity and velocity, respectively, $f(x)$ and $w(x)$ are deterministic forcing term and boundary condition, respectively. 
We denote by $W(D)$ a Banach space and assume that the underlying random input data are properly chosen, such that the corresponding stochastic system \eqref{Lauf} is well-posed and has a unique solution $u(x, \omega) \in L_{\mathbb{P}}^2(\Omega; W(D))$, where the function space
\begin{equation*}
\begin{aligned}
 L_{\mathbb{P}}^2(\Omega; W(D))
  :=\bigg\{ &u:   \overline{D} \times \Omega \rightarrow \mathbb{R} \; \Big|\;  u \mbox{ is strongly measurable}\\
  & \hspace{0cm}\mbox{ and }  \int_{\Omega}\|u\|^2_{W(D)} \,d\mathbb{P}(\omega) < + \infty \bigg\},
\end{aligned}
\end{equation*}
consists of Banach-space valued functions that have finite second moments. Two examples posed in this setting are given below:
\begin{example}[The diffusion equation with random diffusivity]\label{ex1}
{Find a function $u: \overline{D} \times \Omega\rightarrow \mathbb{R}$, such that it holds $\mathbb{P}$-a.e.~in $\Omega$
\begin{equation}\label{eq:ellip}
\left\{
\begin{aligned}
- \nabla \cdot (a(x, \omega)\nabla u(x, \omega))
& =\, f(x)\;\;\;
\mbox{  in  }\; D, \\
u(x, \omega) 
&= \,0\;\;\;\;\;\;\;\;
\mbox{  on  }\;  \partial D,\\
\end{aligned}\right.
\end{equation}
 where the well-posedness is guaranteed in $ L_{\mathbb{P}}^2(\Omega; W(D))=L_{\mathbb{P}}^2(\Omega; H_0^1(D))$ with $f(x) \in L^2(D)$ and $a(x, \omega)$ uniformly elliptic, i.e., for $\mathbb{P}$-a.e. $\omega\in\Omega$,
\begin{equation}\label{eq:uniformellipticity}
	a_{\min} \leq \|a(x,\omega)\|_{L^{\infty}(D)}  \,\,\,\mbox{with}\,\,\, a_{\min} \in (0,\infty).
\end{equation}
}
\end{example}

\begin{example}[The convection-dominated transport with random velocity]\label{ex2}
{Find a function $u: \overline{D} \times \Omega\rightarrow \mathbb{R}$, such that it holds $\mathbb{P}$-a.e.~in $\Omega$
\begin{equation}\label{eq:conv}
\left\{
\begin{aligned}
-  \varepsilon \Delta u(x, \omega) +  \bm{b}(x, \omega) \cdot \nabla u (x, \omega)
& =\, f(x)\;\;\;\;\;\;\;\;
\mbox{  in  }\; D, \\
u(x, \omega) 
&= \,w(x)\;\;\;
\mbox{  on  }\;  \partial D,\\
\end{aligned}\right.
\end{equation}
where $\varepsilon >0$ and the boundary condition $w(x)$ is defined by
\[
w(x) := \left\{
\begin{aligned}
& 1, \quad x \in \mathscr{D},\\
& 0, \quad x \in \partial D \backslash \mathscr{D},\\
\end{aligned}
\right.
\]
where $\mathscr{D}$ is a subset of the boundary $\partial D$. When $\varepsilon$ is very small, the solution will have a sharp transition layer whose location is determined by the velocity field.

}
\end{example}
The two examples exhibit different aspects of the parametric dependence of the solution $u$ on the random coefficients. 
In Example \ref{ex1}, the random diffusion operator leads to a very smooth solution in both the physical domain $D$ and the parametric domain $\Omega$; in Example \ref{ex2}, the random velocity field $\bm{b}(x, \omega)$ will result in sharp transitions of the solution $u$ in the domain $D$, and such irregular behavior will propagate to the stochastic domain $\Omega$. In the next subsection, we introduce the random fields of interest to this effort.

\subsection{The random fields of interest}\label{sec:RF}
We are interested a generic stochastic process on $(\Omega, \mathcal{F}, \mathbb{P})$, denoted by 
\begin{equation}
\eta (x, \omega) :   D \times \Omega \rightarrow \mathbb{R}.
\end{equation}
For a fixed $x \in D$, $\eta(x, \cdot)$ is a real-value square integrable random variable, i.e., 
\begin{equation}
\eta (x, \cdot) \in L^2(\Omega,\mathcal{F},\mathbb{P}) := \left\{ X: \Omega \rightarrow \mathbb{R}\; \bigg| \int_\Omega |X(\omega)|^2 d\mathbb{P}(\omega) < \infty \right\},
\end{equation}
where $L^2(\Omega,\mathcal{F},\mathbb{P}) $ is equipped with the inner product $\langle X, Y \rangle_{\mathbb{P}} := \mathbb{E}[XY]$ and the norm $\|X\|_{\mathbb{P}} = \langle X, Y \rangle_{\mathbb{P}}^{1/2}$.
For notational simplicity, we assume that $\mathbb{E}[\eta(x,\cdot)] = 0 $ for all $x \in D$. The covariance function, denoted by
\begin{equation}\label{cov}
\kappa(x,x') := \mathbb{E}[\eta(x,\omega) \eta(x',\omega)],
\end{equation}
is symmetric and bounded as $\eta \in L^2(\Omega, \mathcal{F},\mathbb{P})$.  The random fields $a(x, \omega)$ and/or $\bm{b}(x, \omega)$ in \eqref{Lauf} are (nonlinear) functions of $\eta(x, \omega)$. For example, the diffusion coefficient could be defined as $a(x, \omega) := \exp(\eta(x,\omega))$ to satisfy the assumption in \eqref{eq:uniformellipticity}; and the random velocity (for $d = 2$) could be defined as 
$b_1(x, \omega) := \cos( \eta(x, \omega))$ and $b_2(x, \omega) := \sin(\eta(x, \omega))$. Of particular interest to this effort are the colored noise discussed in Section \ref{color} and the discrete white noise discussed in Section \ref{white}.

\subsubsection{The colored noise}\label{color} When the covariance $\kappa(x,x')$ is continuous in $D \times D$, we define a compact positive self-adjoint operator $K: L^2(D) \rightarrow L^2(D)$, i.e.,
$
K[v](x) := \int_D \kappa(x,x') v(x') dx',
$
such that $K[\cdot]$ has a complete set of eigenvectors $\{\eta_n(x), n \in \mathbb{N}^+\}$ in $L^2(D)$ and real eigenvalues $\{\lambda_n, n \in \mathbb{N}^+\}$. On the other hand, continuity of $\kappa(x,x')$ implies that the random field $\eta(x,\omega)$ is a mean square continuous stochastic process, i.e., $\lim_{\varepsilon \rightarrow 0} \mathbb{E}[(\eta(x+\varepsilon, \omega) - \eta(x, \omega))^2] = 0$. Thus, $\eta(x,\omega)$ can be represented in the space $span\{ \xi_n(x), n \in \mathbb{N}^+\}$ as
\begin{equation}\label{e10}
\eta(x,\omega) = \sum_{n = 1}^\infty \sqrt{\lambda_n}\, \xi_n(x) y_n(\omega),
\end{equation}
where the random variables $y_n$ are defined by
\[
y_n(\omega) := \int_D \eta(x,\omega) \xi_n(x) dx \;\;\; \text{for} \;\;\; n = 1, 2, \ldots, 
\]
satisfying
$\mathbb{E}[y_n] = 0$, $\mathbb{E}[y_n y_{m}] = \delta_{nm}$ and $\mathbb{V}ar[y_n] = 1$. 
Note that, as long as the correlation is not zero, the eigenvalues will decrease with $n$. The decay rate depends on the covariance function. 
We can approximate $\eta(x,\omega)$ by truncating its Karhunen-Lo\`{e}ve (KL) expansion of the form
\begin{equation}\label{e3}
\eta_N(x, \omega) :=  \sum_{n = 1}^N \sqrt{\lambda_n}\, \xi_n(x) y_n(\omega),
\end{equation}
such that $\eta(x,\omega) \approx \eta_N(x, \omega)$ can be approximately simulated by drawing samples of the $N$-dimensional random vector $\bm y := (y_1, \ldots, y_N)^{\top}$. The representation in \eqref{e3} can be viewed as an approximate parameterization of the original random field $\eta(x, \omega)$. For convenience, we can write the truncated KL expansion as a function of $\bm y$, i.e., $\eta_N(x,\bm y)$. Figure \ref{fig20} shows three snapshots of the random field in \eqref{aaa} with covariance function in \eqref{cor_len} and $L = 0.25$.

Even though there exists a KL expansion for any mean-square continuous stochastic process, it is not easy to obtain the joint probability 
distribution of the random vector $\bm y$. In practice, Gaussian random fields are the most widely used models, which assume that $\eta(x,\omega)$ is a Gaussian random variable for any fixed $x \in D$. In this case, the non-correlation property $\mathbb{E}[y_ny_m] = \delta_{nm}$ leads to {\em independence} of $y_1, \ldots, y_N$, thus $\eta_N(x, \bm y)$ can be easily simulated by sampling the $N$-dimensional standard Gaussian distribution.

\subsubsection{The discrete white noise}\label{white}
We assume the spatial domain $D$ is the union of non-overlapping sub-domains $D_n^{\rm WN}$ for $n = 1,\ldots,N$, i.e.,
\begin{equation}\label{wn}
D = \bigcup_{n=1}^N \overline{D_n^{\rm WN}} \; \; \text{ and }\;\; D_n^{\rm WN} \cap D_{m}^{\rm WN}  = \emptyset \;\; \text{ if }\; m \not= n.
\end{equation}
Then, $\eta(x,\omega)$ is defined as a random variable in each sub-domain $D_n^{\rm WN}$, i.e.,
\begin{equation}\label{e4}
 \eta_N(x, \omega) := \sum_{n=1}^{N}  \mathbf{1}_{D_n^{\rm WN}}(x) y_{n}(\omega),
\end{equation}
where $\mathbf{1}_{D^{\rm WN}_n}(x)$ denotes the indicator function of the sub-domain $D_n^{\rm WN}$. The random variables $y_1,\ldots, y_N$ could be either correlated or independent, bounded or unbounded. Once the joint probability distribution of $\bm y$ is defined, realizations of $\eta(x, \omega) = \eta_N(x, \omega)$ can be generated by drawing samples using Monte Carlo methods. Figure \ref{fig44} shows three snapshots of the discrete white noise with $N = 256$, $\mathbb{E}[y_n] = 0$ and $\mathbb{V}ar[y_n] = 1$.


\section{The deterministic domain decomposition method}\label{sec3}
We briefly review the deterministic domain decomposition method for the PDE in \eqref{Lauf} for a fixed parameter $\omega \in \Omega$, as well as set up necessary notations for discussing our method in Section \ref{sec:DD_MR}. 
%
We define a triangle mesh for the domain $D$, denoted by $\mathcal{T}_h$, which satisfies regular geometric conditions. In this work, we utilize a finite element space, denoted by $X_h\subset H^1(D)$, consisting of piecewise linear continuous basis functions on the conforming triangles of $\mathcal{T}_h$. Note that we use $J$ to represent the total degrees of freedom of 
$X_h$. We denote by $X_h^{0}$ the homogenous counterpart of $X_h$.

For a fixed $\omega \in \Omega$, a finite element scheme for the PDE in \eqref{Lauf} is described as: seek a function $u_h\in X_h$ satisfying $u_h|_{\partial D}=w(x)$ and
\begin{eqnarray}\label{eq:FEM}
\mathcal{A}(u_h,\nu;\omega)=(f,\nu),\;\; \forall \nu\in X_h^0,
\end{eqnarray}
where $\mathcal{A}(u_h, \nu; \omega)$ is a parameterized bilinear form. Then, the solution $u_h$ can be represented in $X_h$ in the form of
%
%
 \[
 u_h(x,\omega)=\sum_{j=1}^J U_j(\omega)\,\psi_j(x),
 \]
where $\left\{\psi_j(x)\right\}_{j=1}^{J}$ is a basis of $X_h$ and ${\mathbf U}(\omega):=\left(U_1(\omega),\cdots,U_J(\omega)\right)^\top$ is the vector of nodal values. It should be noted that different PDEs may need different definitions of the bilinear form in \eqref{eq:FEM}. For example, for the diffusion problem in Example \ref{ex1}, the bilinear form is simply $\mathcal{A}(\mu,\nu;\omega) := (a(x,{\omega})\nabla \mu,\nabla \nu)$; for the convection-dominated transport problem in Example \ref{ex2}, we might need to use the streamline-upwind petrov-Galerkin (SUPG) method, i.e.,
\begin{eqnarray*}
\mathcal{A}(\mu,\nu;\omega):=(\varepsilon\nabla \mu,\nabla \nu)+({\bm b}\cdot\nabla \mu,\nu)+ \sum_{\tau \in \mathcal{T}_h} \delta_\tau(-\varepsilon\Delta \mu +  {\bm b}\cdot\nabla \mu - f,{\bm b}\cdot\nabla \nu)_\tau,
\end{eqnarray*}
to stabilize the finite element scheme, where $\delta_\tau$ is a nonnegative stabilization parameter, and $(\cdot, \cdot)_{\tau}$ is the inner product within the triangle $\tau$.


Now we introduce the domain decomposition. We decompose the physical domain $D$ into $S$ non-overlapping sub-domains, denoted by $D_s^H, s=1,2,\dots,S$, such that
\[
D = \bigcup_{s=1}^N \overline{D_s^{H}} \; \; \text{ and }\;\; D_s^{H} \cap D_{t}^{H}  = \emptyset \;\; \text{ if }\; s \not= t,
\]
and we denote the collection of all the edges and interfaces by
\[
\mathcal{E}^{H} := \bigcup_{s=1}^S \partial D_s^H \backslash \partial D.
\]
One example of such decomposition can be found in Figure \ref{fig:DH-basis}(a), where $D$ is a 2-dimensional square domain. Nevertheless, our method can be used for domains with more complicated geometries, as long as the decomposition is embedded in the triangle mesh.
 %
Based on the embeddedness , we
restrict $\mathcal{T}_{h}$ and $X_h$ in $\overline{D_s^H}$, and define
\[
\mathcal{T}_{h,s}:=\mathcal{T}_h\cap \overline{D_s^H}\quad  \text{ and } \;\; X_{h,s}:=X_h|_{\overline{D_s^H}},
\]
where $J_s$ is the degrees of freedom of ${X}_{h,s}$.

Within each sub-domain, we can write out a local weak formulation
$
\mathcal{A}(u_h,\nu;\omega)_s=(f,\nu)_s, \,\forall \nu\in X_{h,s},
$
which immediately leads to a local algebraic equation
\begin{equation}\label{locfe}
\mathbb{A}_s {\mathbf U}_s = {\mathbf f}_s,
\end{equation}
where $\mathbb{A}_s$ is the local stiffness matrix, ${\bf f}_s$ is the local right-hand side vector, and $\mathbf{U}_s := (U_{s,1}, \ldots, U_{s,J_s})^{\top}$ is the vector of local nodal values in $\overline{D_s^H}$. The system in \eqref{locfe} is singular due to the lack of a boundary condition. Thus, we divide the components in $\mathbf{U}_s$ into two groups, i.e.,
\[
\mathbf{U}_s := 
\begin{pmatrix}
\mathbf{U}_s^0 \vspace{0.2cm}\\ 
\mathbf{U}_s^{\rm b}
\end{pmatrix},
\]
where $\mathbf{U}_s^0$ and $\mathbf{U}_s^{\rm b}$ are the nodal values in the interior and on the boundary of $D_s^H$, respectively. \footnote{Note that the superscript ${}^{\rm 0}$ of a vector or a matrix indicates that the entries associate with the nodal values in the interior of a sub-domain. Analogously, the superscript ${}^{\rm b}$ indicates the association with nodal values on interfaces. 
}
Then, we can recast the algebraic equation \eqref{locfe} in the following form
\begin{eqnarray}\label{AF}
\begin{pmatrix}
\mathbb{A}_s^{00} & \mathbb{A}_s^{0 \rm b}\vspace{0.2cm}\\
\mathbb{A}_s^{\rm b 0}& \mathbb{A}_s^{\rm bb}
\end{pmatrix} 
\begin{pmatrix}
{\bf U}_s^0\vspace{0.2cm}\\
{\bf U}_s^{\rm b}
\end{pmatrix}
=
\begin{pmatrix}
{\bf f}_s^0\vspace{0.2cm}\\
{\bf f}_s^{\rm b}
\end{pmatrix},
\end{eqnarray}
where ${\bf f}_s^0$ and ${\bf f}_s^{\rm b}$ are the right hand vectors corresponding to unknown $\mathbf{U}_s^0$ and $\mathbf{U}_s^{\rm b}$. 

The system in \eqref{AF} can be further manipulated to eliminate the interior unknowns $\mathbf{U}_s^0$ by representing them using $\mathbf{U}_s^{\rm b}$, i.e.,
\begin{equation}\label{schur}
\mathbf{U}_s^0 = (\mathbb{A}_s^{00} )^{-1} \left({\bf f}_s^0 - \mathbb{A}_s^{0 \rm b} \mathbf{U}_s^{\rm b}\right),
\end{equation}
as well as define the local Schur complement
\begin{equation}\label{Bg}
\mathbb{B}_s := \mathbb{A}_s^{\rm bb} - \mathbb{A}_s^{\rm b 0} \left(\mathbb{A}_s^{00}\right)^{-1} \mathbb{A}_s^{0 \rm b}
 \;\; \text{ and } \;\;
{\bf g}_s := \mathbf{f}_s^{\rm b} - \mathbb{A}_s^{\rm b 0} \left(\mathbb{A}_s^{00}\right)^{-1}  \mathbf{f}_s^0.
\end{equation}
Substituting the relation \eqref{schur} into the linear system \eqref{locfe}, we can assemble a global system to solve the unknowns on the interfaces. Specifically, we need to define a manipulation matrix 
$\mathbb{T}_s$ for each sub-domain, and assemble 
\begin{equation}
\mathbb{B} := \sum_{s=1}^S \mathbb{T}_s^{\top} \mathbb{B}_s\mathbb{T}_s, \quad \mathbf{g} := \sum_{s=1}^S \mathbb{T}_s^{\top} \mathbf{g}_s, \quad \mathbf{U}^{\rm b} :=  \sum_{s=1}^S \mathbb{T}_s^{\top} \mathbf{U}_s^{\rm b},
\end{equation}
where the matrix $\mathbb{T}_s$ is used to put the entries of $\mathbb{B}_s$, $\mathbf{g}_s$ and $\mathbf{U}_s^{\rm b}$ to the correct locations in $\mathbb{B}$, $\mathbf{g}$ and $\mathbf{U}^{\rm b}$, respectively. Note that the size of the square matrix $\mathbb{B}$ is smaller than $\sum_{s=1}^S {\rm dim}(\mathbf{U}_s^{\rm b})$, due to shared interfaces between sub-domains. 
After this, $\mathbf{U}^{\rm b}$ can be obtained by solving the condensed system
\begin{equation}\label{globalsys}
\mathbb{B} \mathbf{U}^{\rm b} = \mathbf{g},
\end{equation}
and $\mathbf{U}_s^0$ can be recovered by substituting $\mathbf{U}_s^{\rm b}$ into \eqref{schur}.

Our goal is to reduce the DD method in the stochastic setting. It is easy to see that the main cost of assembling the condensed matrix $\mathbb{B}_s$ in \eqref{Bg} lies in the inversion of $\mathbb{A}_s^{00}$, especially on a very fine triangle mesh. We either need to compute the real inverse of $\mathbb{A}_s^{00}$, or solving the linear system $\mathbb{A}_s^{00}\mathbf{v} = \mathbf{w}$ with $J_s$ different right-hand sides (see \cite{Toselli:2005hi} for details), both of which are very time-consuming. In the stochastic setting, entries of $\mathbb{B}_s$ are functions of the random parameters, so that the inefficient computation of 
$\mathbb{B}_s$ for a large number of parameter samples is the bottleneck of applying the DD strategy to the parametric PDEs. Thus, how to efficiently approximate $\mathbb{B}_s$ in the parameter space is the focus of the next section.

\section{The domain-decomposition model reduction method}\label{sec:DD_MR}
We will describe the details of the proposed DDMR method in this section. The decomposition of random fields will be discussed in Section \ref{sec:rf}; the offline and the online procedures will be discussed in Section \ref{offline} and \ref{online}, respectively.

\subsection{Decomposition of the random fields}\label{sec:rf}
We intend to decompose the random fields of interest into the following form
\begin{equation}\label{locRF}
\eta(x, \omega) \approx \eta_{\bm N}(x, \omega) := \sum_{s =1}^S \eta^{\rm loc}_{s, N_s}(x, \omega)   \mathbf{1}_{D_s^H}(x),
\end{equation}
where $\bm N := (N_1, \ldots, N_S)$ is the vector of dimensions of local random fields, $  \mathbf{1}_{D_s^H}(x)$ is the indicator of the sub-domain $D_s^H$, and $\eta^{\rm loc}_{s, N_s}(x, \omega)$ is the local random field with the support $D_s^H$. Both types of the random fields introduced in Section \ref{sec:RF} can be decomposed and approximated by the form in \eqref{locRF}. We will discuss the colored noise case in Section \ref{sec:locKL} and the discrete white noise case in Section \ref{sec:locWN}.
  

\subsubsection{Local KL expansion for the colored noise}\label{sec:locKL}
For any mean-square continuous random field $\eta$ defined in $D$, we can restrict it to each sub-domain $D_s^{H} \subset D$ and define 
a local KL expansion using the same covariance function $\kappa$, i.e., 
\begin{equation}\label{e9}
\eta^{\rm loc}_s(x,\omega)  := \sum_{n = 1}^\infty \sqrt{\lambda_{s,n}}\, \xi_{s,n}(x) y_{s,n}(\omega)\; \; \text{ for} \; s = 1, \ldots, S,
\end{equation}
such that $\eta(x,\omega) = \eta_s^{\rm loc}(x,\omega)$ for $x \in D_s^{H}$. Similarly, we can truncate the local KL expansion of $\eta_s^{\rm loc}$ and define the following approximation:
\begin{equation}\label{e11}
\eta^{\rm loc}_{s,N_s}(x,\omega)  := \sum_{n = 1}^{N_s} \sqrt{\lambda_{s,n}}\, \xi_{s,n}(x)\, y_{s,n}(\omega)\; \;\;\forall x \in D_s^{H},
\end{equation}
for $s = 1, \ldots, S$. Analogously, we also write $\eta^{\rm loc}_{s,N_s}(x, \bm y_s)$ as a function of the local random variables $\bm y_s :=(y_{s,1}, \ldots, y_{s,N_s})^{\top}$. 
It is easy to see the restriction of each realization of \eqref{e10} in $D^{H}_s$ corresponds a unique realization of the local representation $\eta_s^{\rm loc}$ in \eqref{e9}, but it is not true for $\eta_N(x,\bm y)$ and $\eta^{\rm loc}_{s,N_s}(x,\bm y_s)$ due to the truncations. In this work, for each sample of $\bm y$ in \eqref{e3}, we would like to find a sample of $\bm y_s$ to minimize the error $\eta_N(x, \bm y) - \eta^{\rm loc}_{s,N_s}(x,\bm y_s)$ in the $L^2$ sense, i.e., solving the optimization problem
\begin{equation}\label{e12}
\bm y_s = \argmin_{\bm v \in \mathbb{R}^{N_s}} \left\|\eta_N(\cdot, \bm y) - \eta^{\rm loc}_{s,N_s}(\cdot,\bm v)\right\|^2_{L^2(D_{s}^{H})}.
\end{equation}
An immediate question about the problem \eqref{e12} is how big the minimized $L^2$ error is. In fact, for any fixed $\bm y$, there exists a $\omega^{*} \in \Omega$ such that $\eta(x, \omega^*) = \eta_s^{\rm loc}(x, \omega^*)$ for any $x \in D_s^{H}$. We denoted by $\bm y_s^{*}(\omega^*)$ the image of $\omega^*$ and substitute $\bm y_s^*(\omega^*)$ into \eqref{e11}. Then, the error minimized by $\bm y_s$ in \eqref{e12} can be estimated by
\[
\begin{aligned}
& \mathbb{E}\left[\left\|\eta_N(\cdot, \bm y) - \eta^{\rm loc}_{s,N_s}(\cdot,\bm y_s)\right\|^2_{L^2(D_{s}^{H})}\right]\\
 \le\;& \mathbb{E}\left[\left\|\eta_N(\cdot, \bm y) - \eta^{\rm loc}_{s,N_s}(\cdot,\bm y_s^*)\right\|^2_{L^2(D_{s}^{H})}\right]\\
 \le\; & \mathbb{E}\left[\left\|\eta_N(\cdot, \bm y(\omega^*)) - \eta(\cdot,\omega^*) \right\|^2_{L^2(D_{s}^{H})}\right]+ \mathbb{E}\left[ \left\|\eta_s^{\rm loc}(\cdot, \omega^*) - \eta^{\rm loc}_{s,N_s}(\cdot,\bm y_s^*(\omega^*))\right\|^2_{L^2(D_{s}^{H})}\right]\\
\le \; &\mathbb{E}\left[ \sum_{n = N+1}^{\infty} \lambda_n\, y_n^2(\omega^*) \right]+ \mathbb{E}\left[ \sum_{n = N_s+1}^{\infty} \lambda_{s,n}\, y_{s,n}^2(\omega^*) \right]\\
= \; & \sum_{n = N+1}^{\infty} \lambda_n + \sum_{n = N_s+1}^{\infty} \lambda_{s,n},
\end{aligned}
\]
which implies that $\eta_{s,N_s}^{\rm loc}(x,\bm y_s)$ will provide a good approximation to $\eta_{N}(x,\bm y)$ for sufficiently large $N$ and $N_s$.

In practice, the first step to solve \eqref{e12} is to obtain the eigenvalues and eigenvectors of the covariance function in $D$ and $D_s^{H}$, respectively. It can be done analytically for certain type of covariance functions, e.g., exponential and Gaussian, or numerically using Galerkin projection (see \S 2.1 in \cite{2010smuq.book.....L}) and efficient solvers for eigenvalue problems, e.g., ARPACK\footnote{http://www.caam.rice.edu/software/ARPACK/}. After that the problem in \eqref{e12} can be implemented using discrete least squares method. Specifically, we can draw a set of $T$ uniformly distributed random samples, denoted by $\{x_{s,i}\}_{i = 1}^T$, in the sub-domain $D_s^{H}$, and formulate the following discrete least squares (DLS) problem
\begin{equation}\label{e13}
\bm y_s = \argmin_{\bm v \in \mathbb{R}^{N_s}} \sum_{i = 1}^T\left|\eta_N(x_{s,i}, \bm y) - \eta^{\rm loc}_{s,N_s}(x_{s,i},\bm v)\right|^2;
\end{equation}
the optimal choice of $\bm y_s$ can be obtained by solving the normal system 
\begin{equation}\label{rfDLS}
(\mathbf \Xi_s^{\top} \mathbf \Xi_s) \bm y_s = \mathbf \Xi_s^{\top} \bm \eta_{N},
\end{equation}
where
\[
\begin{aligned}
& [\mathbf \Xi_s]_{ij} := \sqrt{\lambda_{s,j}}\, \xi_{s,j}(x_{s,i}),\;\; i = 1, \ldots, T \text{ and } j = 1, \ldots, N_s,\\
& \bm \eta_N := \left(\eta_N(x_{s,i},\bm y), \ldots, \eta_N(x_{s,T},\bm y)\right)^{\top},
\end{aligned}
\]
are a $T\times N_s$ matrix and a $T$-dimensional vector, respectively. Note that the number of samples $T$ needs to be bigger than $N_s$ to guarantee the stability of the DLS method. After solving the DLS problems in all sub-domains, we can construct an approximate global KL expansion by substituting all the local expansions in \eqref{e11} into \eqref{locRF}. 

\begin{remark}
Since the colored noise is discretized on the triangle mesh $\mathcal{T}_h$ as a piecewise constant function, we do not need to force continuity of the local KL expansion obtained by solving \eqref{rfDLS}, as long as the domain decomposition is embedded in the triangle mesh $\mathcal{T}_h$.
\end{remark}

\begin{remark}
Since this effort does not focus on improving the KL expansion generation algorithms, the strategy used here is feasible but surely not optimal. It becomes inefficient when the triangle mesh $\mathcal{T}_h$ becomes very dense. In that case, the cost of decomposing the covariance function, i.e., solving a Fredholm integral equation of the second kind, is very time consuming. In fact, a DD-based parallel KL expansion generator was developed in \cite{Contreras:H77YYtfe}, which can be directly applied in our setting to improve the efficiency of this step. 
\end{remark}

\subsubsection{Decomposition of the discrete white noise}\label{sec:locWN}
A straightforward way to decompose the discrete white noise in \eqref{wn} is to align both decompositions, i.e., letting $S =N$ and $D_n^{\rm WN} = D_s^{H}$ for $s=n = 1, \ldots, S$. In this case, each sub-domain only involves one random parameter. This is the strategy we will use for the numerical examples, as we only consider the case that $N$ in \eqref{e4} is finite. In this case, the operator of each sub-domain problem will depend on one random parameter, which successfully avoid the curse of dimensionality.

\subsection{The offline procedure}\label{offline}
The purpose of the offline procedure is to construct all the components that are needed in the reduced model, as well as finish all the expensive computation, i.e., any operation whose complexity increases as the triangle mesh size $h$ decreases. Details about our offline procedure are given in Section \ref{sec:data}-\ref{sec:DLS} and a short summary is given in Section \ref{offsum}.

\subsubsection{Training data generation}\label{sec:data} We need to generate two sets of training data, one for solution, and another one for the local stiffness matrices. 
We take the colored noise case as an example to describe how to generate training data. The described procedure can be directly applied in the case of having the discrete white noise. 

To generate data of $u_h$, we set the dimension $N$ in \eqref{e3} sufficiently large, such that the error $\eta-\eta_N$ can be neglected. According to the definition of $\eta_N$ in \eqref{e3}, we sample the multi-variate Gaussian distribution $\mathcal{N}(0, \mathbb{I})$ to generate $K_u$ realizations of $\eta_N$, denoted by
$
\{ \eta_N(x, \bm y(\omega_k)),  k = 1, \ldots, K_u\},
$
each of which is stored as a piecewise constant function on the mesh $\mathcal{T}_h$. Then, we substitute $\{ \eta_N(x, \bm y(\omega_k)),  k = 1, \ldots, K_u\}$ into the weak formulation in \eqref{eq:FEM}, and compute the set of realizations of the PDE solutions, denoted by
$
\{u_h(x, \omega_k), k = 1, \ldots, K_u\},
$
or equivalently all the unknowns,
\begin{equation}\label{data}
\left\{\mathbf{U}_s(\omega_k), s = 1, \ldots, S, k = 1, \ldots, K_u\right\},
\end{equation}
in \eqref{locfe}. Note that the data in \eqref{data} is obtained by using realizations of the global expansion $\eta_N$ and the finite element formulation in \eqref{eq:FEM}, such that $\mathbf{U}_s(\omega_k)$ does not contain the projection error caused by \eqref{e13}. Note that, this procedure is used to generate not only the training data, but also the validation data to test the performance of our method in the numerical examples in Section \ref{example}. 

Next, we discuss how to generate training data for the local stiffness matrices. To do this, we instead sample from the local random variables $\bm y_s$. It is known that $\bm y_s$ in \eqref{e11} also follow the multivariate standard normal distribution $\mathcal{N}(0, \mathbb{I})$. However, to make use of Legendre polynomials in Section \ref{sec:DLS}, we draw $K_y$ samples of $\bm y_s$ {\em uniformly} in a bounded rectangle domain $\Gamma_s \subset \mathbb{R}^{N_s}$. The domain $\Gamma_s$ will be set large enough such that the probability of having a sample $\bm y_s \sim \mathcal{N}(0, \mathbb{I})$ fall outside $\Gamma_s$ is smaller than a tolerance. Due to the low-dimensionality of $\bm y_s$, the number of samples of $\bm y_s$ that fall in the low probability region of $\mathcal{N}(0, \mathbb{I})$ will not be relatively small. On the other hand, an alternative strategy is to use the results on polynomial approximations in irregular domains in \cite{Adcock:2018ve}.
Once $\{\bm y_s(\omega_k), k = 1, \ldots, K_y\}$ is generated, we substitute them into the local KL expansion in \eqref{e11} to obtain realizations
\begin{equation}\label{Ky}
\Big\{ \eta_{s,N_s}^{\rm loc}(x, \bm y_s(\omega_k)),  k = 1, \ldots, K_y\Big\}.
\end{equation}
Substituting such set into the stiffness matrix in \eqref{AF}, we can assemble an approximate local stiffness matrix, denoted by
\begin{eqnarray}\label{appAF}
 \widetilde{\mathbb{A}}_s(\bm y_s(\omega_k)) :=\begin{pmatrix}
\widetilde{\mathbb{A}}_s^{00} (\bm y_s(\omega_k))& \widetilde{\mathbb{A}}_s^{0 \rm b}(\bm y_s(\omega_k))\vspace{0.2cm}\\
\widetilde{\mathbb{A}}_s^{\rm b 0}(\bm y_s(\omega_k))& \widetilde{\mathbb{A}}_s^{\rm bb}(\bm y_s(\omega_k)) 
\end{pmatrix}\approx {\mathbb{A}}_s(\bm y(\omega_k)).
\end{eqnarray}
Note that the error between $\widetilde{\mathbb{A}}_s(\bm y_s(\omega_k))$ and ${\mathbb{A}}_s(\bm y(\omega_k))$ results from the projection in \eqref{e13} from the truncated global expansion $\eta_{N}$ onto the truncated local expansion $\eta_{s,N_s}^{\rm loc}$. Nevertheless, for the discrete white noise introduced in Section \ref{white}, we have ${\mathbb{A}}_s(\bm y(\omega_k)) = \widetilde{\mathbb{A}}_s(\bm y_s(\omega_k)), k = 1, \ldots, K_y$, because there is no global-local projection in the decomposition.

\subsubsection{Constructing reduced global and local linear systems} 
\begin{figure}[h!]
\centering
\begin{tabular}{cc}
  \includegraphics[width=.33\textwidth]{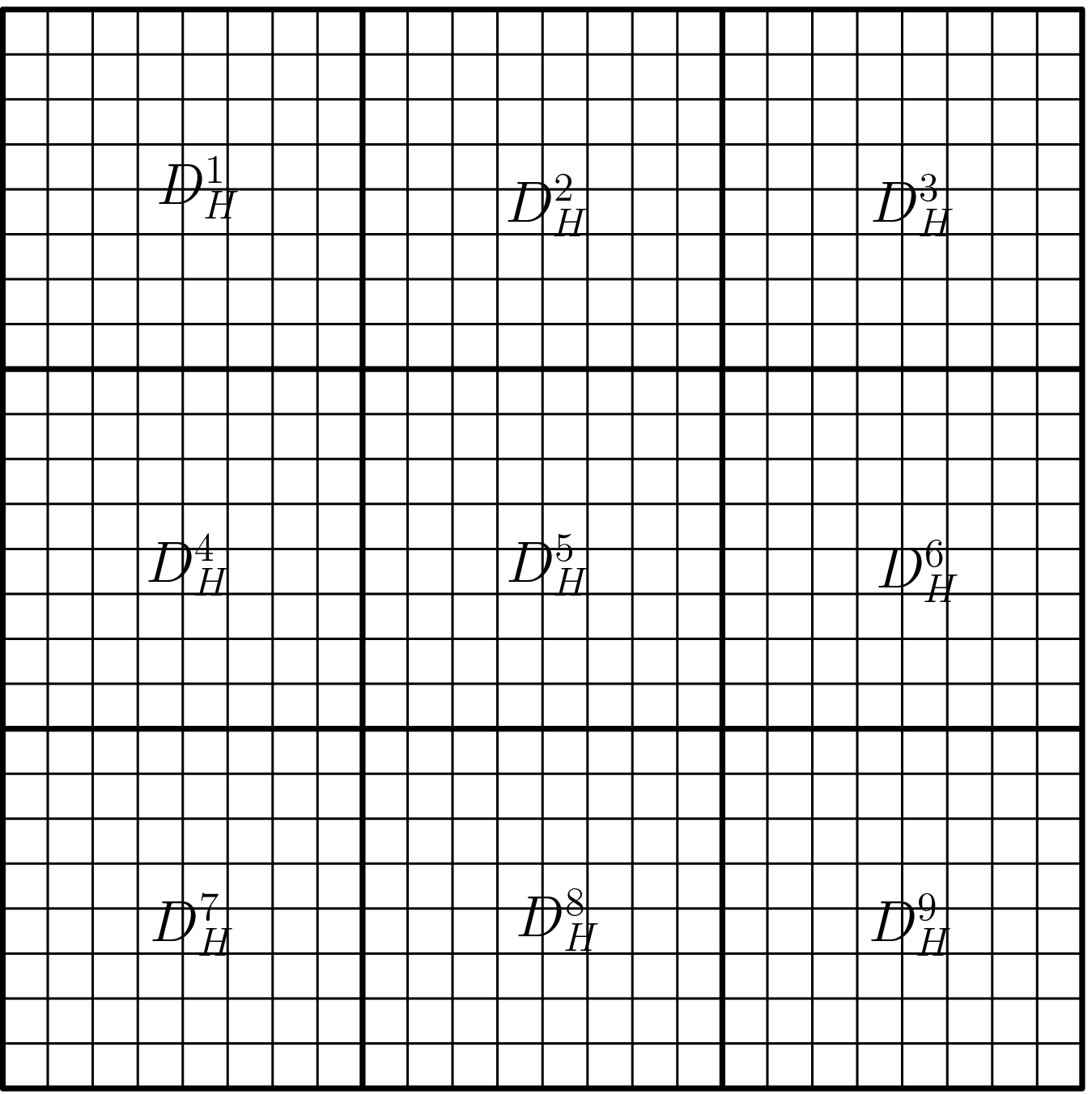}&
   \includegraphics[width=.36\textwidth]{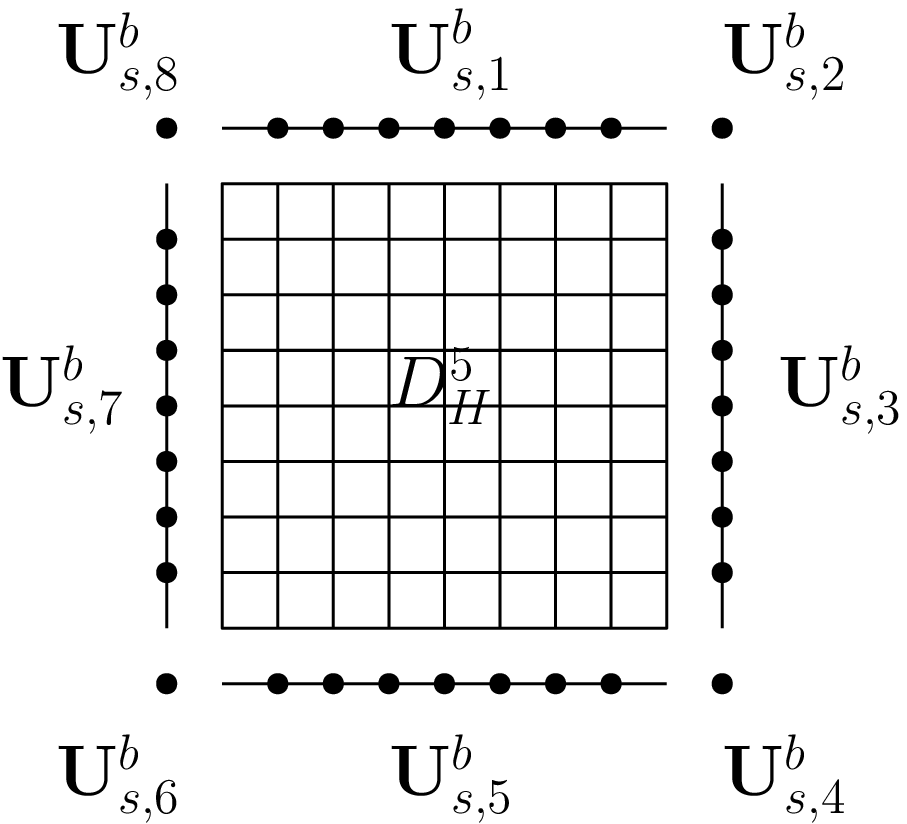} \\
  \footnotesize (a) & \footnotesize(b) 
\end{tabular}
\vspace{-0.2cm}
\caption{(a) Illustration of coarse blocks $D^H_s$ ($s=1,\ldots, 8$); (b) Illustration of how to divide the interface nodal values into different groups. We divide the interface nodal values of $D_5^H$ into 8 non-overlapping groups, i.e., $E_s = 8$, where $\mathbf{U}_{s,1}^{\rm b},\mathbf{U}_{s,3}^{\rm b},\mathbf{U}_{s,5}^{\rm b},\mathbf{U}_{s,7}^{\rm b}$ represent the unknowns on the edges and $\mathbf{U}_{s,2}^{\rm b},\mathbf{U}_{s,4}^{\rm b},\mathbf{U}_{s,6}^{\rm b},\mathbf{U}_{s,8}^{\rm b}$ represent the unknowns at vertices.}\label{fig:DH-basis}
\end{figure}

To reduce the size of the equation in \eqref{globalsys}, we define a set of snapshots for  each interface of each sub-domain, 
\begin{equation}\label{s}
\mathbb{V}_{s,j}^{\rm b} := [\mathbf{U}_{s,j}^{\rm b}(\omega_1), \cdots, \mathbf{U}_{s,j}^{\rm b}(\omega_K)] \;\; \text{ for }\; j = 1, \ldots, E_s, \;\;s = 1, \ldots, S,
\end{equation}
where $E_s$ denotes the number of {\em non-overlapping} groups the nodal values on $\partial D_s^H \backslash \partial D$ are divided into, and $\mathbf{U}_{s,j}^{\rm b}$ denotes the vector of the nodal values of the $j$-th group on $\partial D_s^H \backslash \partial D$.
Taking Figure \ref{fig:DH-basis}(b) as an example, 
we divide the interface nodal values of $D_5^H$ into 8 non-overlapping groups, i.e., $E_s = 8$, where $\mathbf{U}_{s,1}^{\rm b},\mathbf{U}_{s,3}^{\rm b},\mathbf{U}_{s,5}^{\rm b},\mathbf{U}_{s,7}^{\rm b}$ represent the unknowns on the edges and $\mathbf{U}_{s,2}^{\rm b},\mathbf{U}_{s,4}^{\rm b},\mathbf{U}_{s,6}^{\rm b},\mathbf{U}_{s,8}^{\rm b}$ represent the unknowns at vertices. Then, we apply SVD to $\mathbb{V}_{s,j}^{\rm b}$ and generate a reduced basis based on a prescribed threshold, i.e.,
\begin{equation}\label{sing}
\widehat{\mathbb{V}}_{s,j}^{\rm b} := [\mathbf{V}_{s,j,1}^{\rm b}, \cdots, \mathbf{V}_{s,j,M_{s,j}}^{\rm b}],
\end{equation}
where $M_{s,j} < K_u$ is the number of retained left singular vectors. Note that, it is possible that $\mathbb{V}_{s,j}^{\rm b} = \mathbb{V}_{s',j'}^{\rm b}$ when $D_s^H$ and $D_{s'}^H$ share the interface. In this case, we only need to apply SVD to each interface once. 

Once all reduced bases are generated for each sub-domain, we can define the following reduced basis
\begin{equation}\label{vb}
\widehat{\mathbb{V}}_s^{\rm b} := 
\begin{pmatrix}
\widehat{\mathbb{V}}_{s,1}^{\rm b} & & \\
                                & \widehat{\mathbb{V}}_{s,2}^{\rm b} & \\
                                & &  \ddots & \\
                                & & &  \widehat{\mathbb{V}}_{s,E_s}^{\rm b}
\end{pmatrix},
\end{equation}
such that the unknowns $\mathbf{U}_s^{\rm b}$ can be approximated by projecting it onto $\widehat{\mathbb{V}}_s^{\rm b}$, i.e.,
\begin{equation}\label{rb1}
\widehat{\mathbf{U}}_s^{\rm b} = \widehat{\mathbb{V}}_s^{\rm b}\; \widehat{\mathbf{C}}_s^{\rm b},
\end{equation}
where $\widehat{\mathbf{C}}_s^{\rm b}$ is the coefficient vector of size $\sum_{j=1}^{E_s} M_{s,j}$. 

Similarly, we can define a set of snapshots for the interior nodal values of each sub-domain, i.e.,
$
\mathbb{V}_s^0 = [\mathbf{U}_s^0(\omega_1), \cdots, \mathbf{U}_s^0(\omega_{K_u})] $ for $s = 1, \ldots, S,
$
and then apply SVD to obtain a reduced basis for the interior unknowns, i.e.,
\begin{equation}\label{sing1}
\widehat{\mathbb{V}}_{s}^0 := [\mathbf{V}_{s,1}^0, \cdots, \mathbf{V}_{s,M_{s}}^0],
\end{equation}
where $M_{s} < K_u$ is the dimension of the reduced basis. Consequently, we can define an approximation to the interior unknowns $\mathbf{U}_s^0$ by projecting it onto $\widehat{\mathbb{V}}_{s}^0$, i.e.,
\begin{equation}\label{rb2}
\widehat{\mathbf{U}}_s^0 = \widehat{\mathbb{V}}_s^0\; \widehat{\mathbf{C}}_s^0,
\end{equation}
where $\mathbf{C}_s^{\rm b}$ is the coefficient vector of size $M_{s}$. 

Now we can assemble reduced versions of \eqref{schur} and \eqref{globalsys}. For $s = 1, \ldots, S$, we define a reduced system of \eqref{AF} of the form
\begin{equation}\label{rAF}
\begin{aligned}
 &\begin{pmatrix}
\widehat{\mathbb{V}}_{s}^{0} & \\
                                & \widehat{\mathbb{V}}_{s}^{\rm b}
\end{pmatrix}^{\top}
\begin{pmatrix}
\widetilde{\mathbb{A}}_s^{00} & \widetilde{\mathbb{A}}_s^{0 \rm b}\vspace{0.2cm}\\
\widetilde{\mathbb{A}}_s^{\rm b 0}& \widetilde{\mathbb{A}}_s^{\rm bb}
\end{pmatrix}
\begin{pmatrix}
\widehat{\mathbb{V}}_{s}^{0} & \\
                                & \widehat{\mathbb{V}}_{s}^{\rm b}
\end{pmatrix}
=
\begin{pmatrix}
(\widehat{\mathbb{V}}_{s}^{0})^{\top}\widetilde{\mathbb{A}}_s^{00}\,\widehat{\mathbb{V}}_{s}^{0} & (\widehat{\mathbb{V}}_{s}^{0})^{\top}\widetilde{\mathbb{A}}_s^{0 \rm b}\,\widehat{\mathbb{V}}_{s}^{\rm b}\vspace{0.2cm}\\
(\widehat{\mathbb{V}}_{s}^{\rm b})^{\top}\widetilde{\mathbb{A}}_s^{\rm b 0}\,\widehat{\mathbb{V}}_{s}^{0}& (\widehat{\mathbb{V}}_{s}^{\rm b})^{\top}\widetilde{\mathbb{A}}_s^{\rm bb}\,\widehat{\mathbb{V}}_{s}^{\rm b}
\end{pmatrix},\\
& \hspace{3.5cm}\begin{pmatrix}
\widehat{\mathbb{V}}_{s}^{0} & \\
                                & \widehat{\mathbb{V}}_{s}^{\rm b}
\end{pmatrix}^{\top} \begin{pmatrix}
{\bf f}_s^0\vspace{0.2cm}\\
{\bf f}_s^{\rm b}
\end{pmatrix}
=
\begin{pmatrix}
(\widehat{\mathbb{V}}_{s}^{0})^{\top}{\bf f}_s^0\vspace{0.2cm}\\
(\widehat{\mathbb{V}}_{s}^{\rm b})^{\top}{\bf f}_s^{\rm b}
\end{pmatrix}.
\end{aligned}
\end{equation}
The system can be further manipulated to eliminate the interior unknowns $\widehat{\mathbf{C}}_s^0$ by representing them using $\widehat{\mathbf{C}}_s^{\rm b}$, i.e.,
\begin{equation}\label{schur2}
\widehat{\mathbf{C}}_s^0 = 
\left((\widehat{\mathbb{V}}_{s}^{0})^{\top} \widetilde{\mathbb{A}}_s^{00}\, \widehat{\mathbb{V}}_{s}^{0}\right)^{-1}\,
\left((\widehat{\mathbb{V}}_{s}^{0})^{\top} {\bf f}_s^0 -(\widehat{\mathbb{V}}_{s}^{0})^{\top} \widetilde{\mathbb{A}}_s^{0 \rm b}\,
\widehat{\mathbb{V}}_{s}^{\rm b}\,\widehat{\mathbf{C}}_s^{\rm b}\right),
\end{equation}
as well as define reduced forms of the matrix and vector in \eqref{Bg}, i.e.,
\begin{equation}\label{rBg}
\widehat{\mathbb{B}}_s := (\widehat{\mathbb{V}}_{s}^{\rm b})^{\top}\widetilde{\mathbb{A}}_s^{\rm bb}\,\widehat{\mathbb{V}}_{s}^{\rm b} 
- (\widehat{\mathbb{V}}_{s}^{\rm b})^{\top}\widetilde{\mathbb{A}}_s^{\rm b 0}\,\widehat{\mathbb{V}}_{s}^{0}
 \left(  (\widehat{\mathbb{V}}_{s}^{0})^{\top} \widetilde{\mathbb{A}}_s^{00}\, \widehat{\mathbb{V}}_{s}^{0}\right)^{-1} 
  (\widehat{\mathbb{V}}_{s}^{0})^{\top} \widetilde{\mathbb{A}}_s^{0 \rm b}  \widehat{\mathbb{V}}_{s}^{\rm b},
\end{equation}
and
\begin{equation}\label{rg}
\widehat{{\bf g}}_s :=(\widehat{\mathbb{V}}_{s}^{\rm b})^{\top} \mathbf{f}_s^{\rm b} -(\widehat{\mathbb{V}}_{s}^{\rm b})^{\top} \widetilde{\mathbb{A}}_s^{\rm b 0}\, \widehat{\mathbb{V}}_{s}^{0}
 \left((\widehat{\mathbb{V}}_{s}^{0})^{\top}\widetilde{\mathbb{A}}_s^{00}\, \widehat{\mathbb{V}}_{s}^{0}\right)^{-1} (\widehat{\mathbb{V}}_{s}^{0})^{\top} \mathbf{f}_s^0.
\end{equation}

Then, we can cancel out $\widehat{\mathbf{C}}_s^0$ and assemble a global system to solve $\widehat{\mathbf{C}}_s^{\rm b}$. Specifically, we need to define another manipulation matrix 
$\widehat{\mathbb{T}}_s$ for each sub-domain, and assemble 
\begin{equation}
\widehat{\mathbb{B}} := \sum_{s=1}^S \widehat{\mathbb{T}}_s^{\top} \widehat{\mathbb{B}}_s\, \widehat{\mathbb{T}}_s, \quad 
\widehat{\mathbf{g}} := \sum_{s=1}^S \widehat{\mathbb{T}}_s^{\top} \widehat{\mathbf{g}}_s, \quad 
\widehat{\mathbf{C}}^{\rm b} :=  \sum_{s=1}^S \widehat{\mathbb{T}}_s^{\top} \widehat{\mathbf{C}}_s^{\rm b},
\end{equation}
where the matrix $\widehat{\mathbb{T}}_s$
is used to put the entries of $\widehat{\mathbb{B}}_s$, $\widehat{\mathbf{g}}_s$ and $\widehat{\mathbf{C}}_s^{\rm b}$ to the correct locations in $\widehat{\mathbb{B}}$, $\widehat{\mathbf{g}}$ and $\widehat{\mathbf{C}}^{\rm b}$, respectively. 
Note that the size of the reduced square matrix $\widehat{\mathbb{B}}$ is smaller than $\sum_{s=1}^S {\rm dim}(\widehat{\mathbf{C}}_s^{\rm b})$, due to shared interfaces between sub-domains. After this, $\widehat{\mathbf{C}}^{\rm b}$ can be obtained by solving 
\begin{equation}\label{globalsys1}
\widehat{\mathbb{B}}\, \widehat{\mathbf{C}}^{\rm b} = \widehat{\mathbf{g}},
\end{equation}
and $\widehat{\mathbf{C}}_s^0$ can be recovered by substituting $\widehat{\mathbf{C}}_s^{\rm b}$ into \eqref{schur2}.

Since each matrix of snapshots ${\mathbb{V}}_s^0$ or ${\mathbb{V}}_{s,j}^{\rm b}$ only covers a small portion of nodal values, its singular value will decay faster than the case of applying SVD to the matrix of snapshots of all nodal values of the global system. Similar to the local KL expansion, the decay rate of singular value of ${\mathbb{V}}_s^0$ or ${\mathbb{V}}_{s,j}^{\rm b}$ will also depend on the size of the sub-domains. A large number of sub-domains will lead to faster decay of singular values, so that a smaller value for $M_{s,j}$ would be sufficient to achieve the prescribed accuracy. Nevertheless, we do not know the optimal domain decomposition strategy to obtain a system in \eqref{globalsys1} with minimal size.

\begin{remark}
The strategy of applying POD to reduce the dimension on the interfaces have been used to develop other static condensation methods, e.g., \cite{Liao:2015gy,Eftang:2013jp,Huynh:2012cf}. Our contribution lies in the integration of this strategy into our methodology to successfully address the challenges of high-dimensionality and irregular behaviors, especially for the convection-dominated PDEs with random velocities in Example \ref{ex2}. 
\end{remark}

\subsubsection{Sparse approximation of the local stiffness matrices}\label{sec:DLS}
So far we managed to reduce the sizes of global and local systems by generating reduced bases in the sub-domains and on the interfaces. The remaining challenge is the cost of assembling the reduced stiffness matrices $\widehat{\mathbb{B}}_s(\omega)$ and the right-hand side $\widehat{\mathbf{g}}_s(w)$ for a large number of samples of $\omega$. To reduce this cost, we propose to construct sparse polynomial approximations to the entries of the following matrices
\begin{equation}\label{DLStar}
\begin{aligned}
&\widehat{\mathbb{A}}_s^{\rm 00}(\bm y_s) := (\widehat{\mathbb{V}}_{s}^{0})^{\top}\widetilde{\mathbb{A}}_s^{00}(\bm y_s)\,\widehat{\mathbb{V}}_{s}^{0}, \qquad
  \widehat{\mathbb{A}}_s^{\rm 0b}(\bm y_s) := (\widehat{\mathbb{V}}_{s}^{0})^{\top}\widetilde{\mathbb{A}}_s^{0 \rm b}(\bm y_s)\,\widehat{\mathbb{V}}_{s}^{\rm b}, \\
&\widehat{\mathbb{A}}_s^{\rm b0}(\bm y_s) := (\widehat{\mathbb{V}}_{s}^{\rm b})^{\top}\widetilde{\mathbb{A}}_s^{\rm b 0}(\bm y_s)\,\widehat{\mathbb{V}}_{s}^{0},\qquad
  \widehat{\mathbb{A}}_s^{\rm bb}(\bm y_s) :=(\widehat{\mathbb{V}}_{s}^{\rm b})^{\top}\widetilde{\mathbb{A}}_s^{\rm bb}(\bm y_s)\,\widehat{\mathbb{V}}_{s}^{\rm b},
 \end{aligned}
\end{equation}
for $s = 1, \ldots, S$, where $\widetilde{\mathbb{A}}_s^{\rm 00},\widetilde{\mathbb{A}}_s^{\rm 0 b},\widetilde{\mathbb{A}}_s^{\rm b0},\widetilde{\mathbb{A}}_s^{\rm bb}$ are defined in \eqref{appAF}. It should be noted that the matrices in \eqref{DLStar} only depend on the local random vector $\bm y_s$ of dimension $N_s$, so that we can exploit the dimension reduction benefit in the sparse approximation. Specifically, we define a sparse polynomial approximation in the bounded domain $\Gamma_s \subset \mathbb{R}^{N_s}$ using Legendre basis, i.e.,
\begin{equation}\label{DLS}
\begin{aligned}
& \widehat{\mathbb{A}}_s^{\rm 00, LS}(\bm y_s) := \sum_{m = 1}^{M^{\rm LS}_{s,p}} \widehat{\mathbb{K}}_{s,m}^{00} \, L_m (\bm y_s) , \qquad 
   \widehat{\mathbb{A}}_s^{\rm 0b, LS}(\bm y_s) := \sum_{m = 1}^{M^{\rm LS}_{s,p}} \widehat{\mathbb{K}}_{s,m}^{\rm 0b} \, L_m (\bm y_s),\\
   & \widehat{\mathbb{A}}_s^{\rm b0, LS}(\bm y_s) := \sum_{m = 1}^{M^{\rm LS}_{s,p}} \widehat{\mathbb{K}}_{s,m}^{\rm b0} \, L_m (\bm y_s) , \qquad 
   \widehat{\mathbb{A}}_s^{\rm bb, LS}(\bm y_s) := \sum_{m = 1}^{M^{\rm LS}_{s,p}} \widehat{\mathbb{K}}_{s,m}^{\rm bb} \, L_m (\bm y_s),
\end{aligned}
\end{equation}
where $L_m(\bm y_s)$ for $m = 1, \ldots, M_{s,p}^{\rm LS}$ are Legendre polynomials expanding the space $\mathcal{P}_{M^{\rm LS}_{s,p}}(\Gamma_s)$, where the subscript $p$ shows the maximum polynomial order of in the space. The coefficient matrices in \eqref{DLS} are computed using the training data generated in \eqref{appAF}. 

In the colored noise case, since the dimension $N_s$ of the local random vector $\bm y_s$ is much smaller than the global dimension $N$, we need a much smaller cardinality $M_{s,p}^{\rm LS}$ to achieve a prescribed accuracy. Moreover, the local dimension $N_s$ can be further reduced by chopping the domain $D$ into more sub-domains. In this work, we use {\em anisotropic} total degree polynomial spaces \cite{Nobile:2008dr} to define $\mathcal{P}_{M^{\rm LS}_{s,p}}(\Gamma_s)$, where the anisotropy is determined based on the singular value decay of the local KL expansion. More advanced method could be used to further exploit the sparsity, even though the simple anisotropic space is sufficient to illustrate the superior performance of our method. 
In the discrete white noise case, we set $D_n^{\rm WN} = D_s^{H}$, i.e., aligning the interfaces with the partition of the noise, so that the local dimension is $N_s = 1$.

A major difference between our method and the existing work, e.g., \cite{Contreras:aZvF515u}, is that we can approximate the local and global Schur complement matrices $\mathbb{B}_s$ and $\mathbb{B}$ without approximating a large number of local problems. Due to inefficiency of inverting $\mathbb{A}_s^{\rm 00}$ and the lack of boundary condition in \eqref{AF}, the exisitng strategy to construct $\mathbb{B}_s$ is to decompose the local problem in \eqref{AF} into a set of $\sum_{j=1}^{E_j} M_{s,j}$ sub-problems and approximate the parametric dependence for all the sub-problems. Nevertheless, in our method, 
we managed to reduce both the sizes of $\widehat{\mathbb{A}}_s^{\rm 00},\widehat{\mathbb{A}}_s^{\rm 0 b},\widehat{\mathbb{A}}_s^{\rm b0},\widehat{\mathbb{A}}_s^{\rm bb}$ to $M_{s} \times M_{s}$, $M_{s} \times M_{s,j}$, $M_{s,j} \times M_{s}$, $M_{s,j} \times M_{s,j}$, respectively, and the dimension of the parametric dependence to $N_s$. Thus, we can directly approximate all the entries of $\mathbb{B}_s$ efficiently, which provides a significant saving to the offline cost. In addition, our strategy also makes it straightforward to handle non-affine random coefficients without using the classic empirical interpolation.

On the other hand, the avoidance of building sparse approximations to local PDE solutions is particularly beneficial in solving the stochastic convection-dominated problem with sharp transitions (see Section \ref{CD}). The sharp transition is caused by the discontinuous boundary condition, and the transition layer moves for different realizations of the random velocity field. This will result in sharp transitions of the parametric dependence $\bm y(\omega)\rightarrow u_h(x, \omega)$. It is well known that approximating irregular functions is very challenging, especially in high-dimensional spaces (e.g., see \cite{Zhang:2016tn,Jakeman:2011fw}). However, the local stiffness matrices $\widehat{\mathbb{A}}_s^{\rm 00},\widehat{\mathbb{A}}_s^{\rm 0 b},\widehat{\mathbb{A}}_s^{\rm b0},\widehat{\mathbb{A}}_s^{\rm bb}$ are not affected by such irregularity, so that convergence of the sparse approximation in \eqref{DLS} will not slow down.

\subsubsection{Summary of the offline procedure}\label{offsum}
The offline procedure discussed in this section can be summarized in Algorithm 1, in the case of having the colored noise. The algorithm for handling the discrete white noise can be obtained by a slight modification. 
\vspace{0.2cm}
\begin{center}
\begin{tabular}{p{0.94\textwidth}}
\hline
\noalign{\smallskip}
{\bf Algorithm 1}: {\em The offline procedure for the colored noise case}\\
\noalign
{\smallskip}\hline
\noalign{\smallskip}
\vspace{-0.3cm}
\begin{spacing}{1.5}
\begin{algorithmic}[1]\label{algorithm1}
\renewcommand{\algorithmicrequire}{\textbf{Input:}}
\renewcommand{\algorithmicensure}{\textbf{Output:}}
\Require $D, f, w, a, \bm b$ in \eqref{Lauf}, triangulation $\mathcal{T}_h$, covariance $\kappa(x,x')$ in \eqref{cov};
\Ensure $\{\lambda_n, \xi_n\}_{n = 1}^N$, $\{\lambda_{s,n}, \xi_{s,n}\}_{n = 1, s = 1}^{N_s,S}$,
$\{(\widehat{\mathbb{V}}_s^{\rm b})^{\top} \mathbf{f}_s^{\rm b},\, (\widehat{\mathbb{V}}_s^{\rm 0})^{\top} \mathbf{f}_s^{\rm 0}\}_{s=1}^S$,
$\{\widehat{\mathbb{V}}_{s}^{\rm 0}\}_{s=1}^S$, $\{\widehat{\mathbb{V}}_{s}^{\rm b}\}_{s=1}^S$,
$\{\widehat{\mathbb{K}}_{s,m}^{00}\}_{s=1,m=1}^{S,M_{s,p}^{\rm LS}}$, $\{\widehat{\mathbb{K}}_{s,m}^{\rm 0b}\}_{s=1,m=1}^{S,M_{s,p}^{\rm LS}}$, $\{\widehat{\mathbb{K}}_{s,m}^{\rm b0}\}_{s=1,m=1}^{S,M_{s,p}^{\rm LS}}$, 
$\{\widehat{\mathbb{K}}_{s,m}^{\rm bb}\}_{s=1,m=1}^{S,M_{s,p}^{\rm LS}}$;
\State Compute eigenvalues $\lambda_n$ and eigenvectors $\xi_n$ for $n = 1, \ldots, N$ for $\eta_N$ in \eqref{e3};
\State Decompose the domain $D$ into $D_s^H$ for $s = 1, \ldots, S$;
\State Compute eigenvalues $\lambda_{s,n}$ and eigenvectors $\xi_{s,n}$ for $\eta_{s,N_s}^{\rm loc}$ in \eqref{e11}; 
\State Generate training data $\left\{\mathbf{U}_s(\omega_k), s = 1, \ldots, S, k = 1, \ldots, K_u\right\}$ in \eqref{data};
\State Use SVD to construct $\widehat{\mathbb{V}}_{s,j}^{\rm b} $ in \eqref{sing} for $j = 1, \ldots, E_s, s = 1,\ldots, S$;
\State Assemble $\widehat{\mathbb{V}}_s^{\rm b}$ in \eqref{vb} for $s = 1, \ldots, S$;
\State Use SVD to construct $\widehat{\mathbb{V}}_{s}^{\rm 0} $ in \eqref{sing1} for $s = 1,\ldots, S$;
\State Compute and store $(\widehat{\mathbb{V}}_s^{\rm b})^{\top} \mathbf{f}_s^{\rm b}$ and $(\widehat{\mathbb{V}}_s^{\rm 0})^{\top} \mathbf{f}_s^{\rm 0}$ for $s = 1, \ldots, S$;
\State Generate training data $\{ \widetilde{\mathbb{A}}_s(\bm y_s(\omega_k)),\, s = 1, \ldots, S, k = 1, \ldots, K_y\}$ in \eqref{appAF};
\State Construct reduced data $\widehat{\mathbb{A}}_s^{\rm 00}(\bm y_s), \widehat{\mathbb{A}}_s^{\rm 0b}(\bm y_s), \widehat{\mathbb{A}}_s^{\rm b0}(\bm y_s), \widehat{\mathbb{A}}_s^{\rm bb}(\bm y_s)$ using \eqref{DLStar};
\State Solve the coefficient matrices $\widehat{\mathbb{K}}_{s,m}^{00}, \widehat{\mathbb{K}}_{s,m}^{\rm 0b}, \widehat{\mathbb{K}}_{s,m}^{\rm b0}, \widehat{\mathbb{K}}_{s,m}^{\rm bb}$ in \eqref{DLS};
\end{algorithmic}
\vspace{-1.1cm}
\end{spacing}
\\
\hline
\end{tabular}
\end{center}
\vspace{0.2cm}

In terms of number of operations, the dominant cost lies in the generation of $\{\mathbf{U}_s(\omega_k), s =1, \ldots, S, k = 1, \ldots, K_u\}$ and
$\{ \widetilde{\mathbb{A}}_s(\bm y_s(\omega_k)),\, s = 1, \ldots, S, k = 1, \ldots, K_y\}$. Both numbers $K_u$ and $K_y$ could be reduced by dividing $D$ into more sub-domains. The number $K_u$ is an important factor on the quality of the solution manifold coverage by the training data. As the size of each sub-domain becomes smaller, the set $\{\mathbf{U}_s(\omega_k), k = 1, \ldots, K_u\}$ only covers a smaller sub-manifold in a lower dimensional space. As such, given sufficient regularity of the manifold, e.g., the PDE in Example \ref{ex1}, we expected that a smaller value of $K_u$ can provide sufficient coverage to achieve a prescribed accuracy. The number $K_y$ is related to the cardinality $M_{s, p}^{\rm LS}$ of the polynomial space. Thus, it is easy to see that a smaller sub-domain will lead to a better anisotropy of the local KL expansion, such that $M_{s, p}^{\rm LS}$ can be further reduced for a given polynomial order $p$. On the other hand, we would like to emphasize again that the domain decomposition is only used to effective dimension reduction, but the convergence of our method does not require the number of sub-domains goes to infinity.

In terms of storage requirement, the sizes of the matrices for storing bases functions, i.e., $\{\xi_n\}_{n = 1}^N$, $\{\xi_{s,n}\}_{n = 1, s = 1}^{N_s,S}$, $\{\widehat{\mathbb{V}}_{s}^{\rm 0}\}_{s=1}^S$ and $\{\widehat{\mathbb{V}}_{s}^{\rm b}\}_{s=1}^S$, depend on the triangle mesh size $h$, which is unavoidable. A major improvement of this effort is that the sizes of the sparse approximation coefficients, i.e.,   
$\widehat{\mathbb{K}}_{s,m}^{00}$, $\widehat{\mathbb{K}}_{s,m}^{\rm 0b}$, $\widehat{\mathbb{K}}_{s,m}^{\rm b0}$, 
$\widehat{\mathbb{K}}_{s,m}^{\rm bb}$ are of size $M_{s} \times M_{s}$, $M_{s} \times M_{s,j}$, $M_{s,j} \times M_{s}$, $M_{s,j} \times M_{s,j}$, respectively, which are, again, independent of triangle mesh size $h$. Thus, the space required to store those coefficients are on the order of $\mathcal{O}(S M_{s, p}^{\rm LS}(M_s + M_{s,j})^2)$ when we use the same $M_{s}, M_{s,j}$ and $M_{s,p}^{\rm LS}$ for all sub-domains and interfaces.
%

\subsection{The online procedure}\label{online}
The online procedure involves how to use outputs of the offline procedure to approximate the solution with the cost independent of the original triangle mesh size $h$. We summarize the online procedure in Algorithm 2, for the case of having the colored noise. The algorithm for handling the discrete white noise can be obtained by a slight modification. 
 
\vspace{0.2cm}
\begin{center}
\begin{tabular}{p{0.94\textwidth}}
\hline
\noalign{\smallskip}
{\bf Algorithm 2}: {\em The online procedure for the colored noise case}\\
\noalign
{\smallskip}\hline
\noalign{\smallskip}
\vspace{-0.3cm}
\begin{spacing}{1.5}
\begin{algorithmic}[1]\label{algorithm2}
\renewcommand{\algorithmicrequire}{\textbf{Input:}}
\renewcommand{\algorithmicensure}{\textbf{Output:}}
\Require 
$\{\lambda_n, \xi_n\}_{n = 1}^N$, $\{\lambda_{s,n}, \xi_{s,n}\}_{n = 1, s = 1}^{N_s,S}$,
$\{(\widehat{\mathbb{V}}_s^{\rm b})^{\top} \mathbf{f}_s^{\rm b},\, (\widehat{\mathbb{V}}_s^{\rm 0})^{\top} \mathbf{f}_s^{\rm 0}\}_{s=1}^S$,
$\{\widehat{\mathbb{V}}_{s}^{\rm 0}\}_{s=1}^S$, $\{\widehat{\mathbb{V}}_{s}^{\rm b}\}_{s=1}^S$,
$\{\widehat{\mathbb{K}}_{s,m}^{00}\}_{s=1,m=1}^{S,M_{s,p}^{\rm LS}}$, $\{\widehat{\mathbb{K}}_{s,m}^{\rm 0b}\}_{s=1,m=1}^{S,M_{s,p}^{\rm LS}}$, $\{\widehat{\mathbb{K}}_{s,m}^{\rm b0}\}_{s=1,m=1}^{S,M_{s,p}^{\rm LS}}$, 
$\{\widehat{\mathbb{K}}_{s,m}^{\rm bb}\}_{s=1,m=1}^{S,M_{s,p}^{\rm LS}}$;
\Ensure $\widehat{\mathbf{C}}^{\rm b, LS}$, $\{\widehat{\mathbf{C}}^{\rm 0, LS}_s\}_{s = 1}^S$;
\State Generate a sample of $\bm y \in \mathbb{R}^N$;
\State Project $\bm y$ to $\bm y_s$ for $s = 1, \ldots, S$ using \eqref{rfDLS};
\State Substitute $\bm y_s$ into \eqref{DLS} to evaluate 
$\widehat{\mathbb{A}}_s^{\rm 00, LS}$, 
$\widehat{\mathbb{A}}_s^{\rm 0b, LS}$,
$\widehat{\mathbb{A}}_s^{\rm b0, LS}$,
$\widehat{\mathbb{A}}_s^{\rm bb, LS}$;
\State Compute $\big(\widehat{\mathbb{A}}_s^{\rm 00, LS}\big)^{-1}$ for $s = 1, \ldots, S$;
\State Construct $\widehat{\mathbb{B}}_s^{\rm LS}(\bm y_s) := \widehat{\mathbb{A}}_s^{\rm bb, LS} - \widehat{\mathbb{A}}_s^{\rm b0, LS}
\big(\widehat{\mathbb{A}}_s^{\rm 00, LS}\big)^{-1} \widehat{\mathbb{A}}_s^{\rm 0b, LS}$ for $s = 1, \ldots, S$;
\State Construct $\widehat{\mathbf{g}}_s(\bm y_s) := (\widehat{\mathbb{V}}_s^{\rm b})^{\top} \mathbf{f}_s^{\rm b} - \widehat{\mathbb{A}}_s^{\rm b0} \big(\widehat{\mathbb{A}}_s^{00}\big)^{-1} (\widehat{\mathbb{V}}_s^{\rm 0})^{\top} \mathbf{f}_s^{\rm 0}$;
\State Assemble $\widehat{\mathbb{B}}^{\rm LS} := \sum_{s=1}^S (\widehat{\mathbb{T}}_s)^{\top}\, \widehat{\mathbb{B}}_s^{\rm LS}\, \widehat{\mathbb{T}}_s, \quad 
\widehat{\mathbf{g}}^{\rm LS} := \sum_{s=1}^S (\widehat{\mathbb{T}}_s)^{\top} \widehat{\mathbf{g}}_s^{\rm LS}$;
\State Solve the global system $\widehat{\mathbb{B}}^{\rm LS}\, \widehat{\mathbf{C}}^{\rm b, LS} = \widehat{\mathbf g}^{\rm LS}$;
\State Assign $\widehat{\mathbf{C}}^{\rm b, LS}$ to $\widehat{\mathbf{C}}^{\rm b, LS}_s$ for $s = 1, \ldots, S$;
\State Recover the local unknowns $ \widehat{\mathbf{C}}^{\rm 0, LS}_s = \big(\widehat{\mathbb{A}}_s^{\rm 00, LS}\big)^{-1}\left[(\widehat{\mathbb{V}}_s^{\rm 0})^{\top} \mathbf{f}_s^{\rm 0} -  \widehat{\mathbb{A}}_s^{\rm 0b, LS}\,\widehat{\mathbf{C}}^{\rm b, LS}_s\right];
$
\end{algorithmic}
\vspace{-1.1cm}
\end{spacing}
\\
\hline
\end{tabular}
\end{center}
\vspace{0.2cm}

The key feature of the online procedure is that the cost of of solving  
$\widehat{\mathbf{C}}^{\rm b, LS}$, $\{\widehat{\mathbf{C}}^{\rm 0, LS}_s\}_{s = 1}^S$ for each sample $\bm y \in \mathbb{R}^N$ is independent of the triangle mesh size $h$.
First, we can see that the cost of mapping $\bm y$ to $\bm y_s$ in each sub-domain involves solving a linear system of size $T\times T$ in \eqref{rfDLS}, where $T > N_s$ only needs to be big enough to guarantee numerical stability.
Second, since the sizes of $\widehat{\mathbb{A}}_s^{\rm 00, LS}$, 
$\widehat{\mathbb{A}}_s^{\rm 0b, LS}$,
$\widehat{\mathbb{A}}_s^{\rm b0, LS}$ and 
$\widehat{\mathbb{A}}_s^{\rm bb, LS}$ are $M_{s} \times M_{s}$, $M_{s} \times M_{s,j}$, $M_{s,j} \times M_{s}$, $M_{s,j} \times M_{s,j}$, respectively, evaluation of those matrices involves a total of 
$
S(M_s + M_{s,j})^2
$
vector-vector multiplications, where the vectors are of size $M_{s,p}^{\rm LS}$.
A major cost lies in the inversion of $\widehat{\mathbb{A}}_s^{\rm 00, LS}$, which requires $\mathcal{O}(SM_s^3)$ operations. In addition, since the matrix $\widehat{\mathbb{B}}_s^{\rm LS}$ is of size $\sum_{j=1}^{E_s}M_{s,j} \times \sum_{j=1}^{E_s}M_{s,j}$, the size of the global matrix $\widehat{\mathbb{B}}^{\rm LS}$ is smaller than $\sum_{s=1}^S  \sum_{j=1}^{E_s}M_{s,j} \times \sum_{s=1}^S  \sum_{j=1}^{E_s}M_{s,j} $ due to shared nodal values, so that the cost of solving $\widehat{\mathbb{B}}^{\rm LS}\, \widehat{\mathbf{C}}^{\rm b, LS} = \widehat{\mathbf g}^{\rm LS}$ is also independent of the triangle mesh size $h$. 

\subsection{Discussion on approximation errors}
This effort focuses on the development of a new domain decomposition method for the PDEs with random inputs, and rigorous error analysis will be conducted in the future work. Nevertheless, it is not difficult to identify the main error sources from Algorithm 1 and 2. Basically, there are four main error sources, i.e., (i) finite element discretization, (ii) discretization of random fields, (iii) reduced bases representation, (iv) sparse approximation of the local stiffness matrices. The first error can be estimated by following the standard finite element analysis. 
The second source only applies to the colored noise cases, where the error comes from the truncations of the global and local KL expansions, as well as the least-squares projection from the global to the local random variables. According to the theories on KL expansion (e.g., see \cite{Schwab:2006isa}), this error can be controlled by increasing the dimension of the KL expansions, i.e., increasing $N$ and $N_s$ in \eqref{e3} and \eqref{e11}, respectively. The error from the third source is related to the Kolmogorov $n$-width (i.e., the optimal error) for a sub-manifold of the solution, i.e., the manifold of the solution on an interface or in a sub-domain. This error is not easy to analyze due to its correlation with the domain partition. The smaller each sub-domain, the faster the $n$-width decays for a local manifold. Nevertheless, such faster $n$-width decay does not necessarily lead to smaller size of the global system matrix $\widehat{\mathbb{B}}$ in \eqref{globalsys1}, as more sub-domains have to be handled. Thus, an important question to be answered in our future work, is how to partition the domain in an efficient way to minimize the size of the reduced global system in \eqref{globalsys1}. At last, the error of approximating the reduced stiffness matrices depend on the regularity of the entries of those matrices with respect to the local parameters. In fact, since the bilinear form in \eqref{eq:FEM} is a linear or quadratic function of the random coefficients, the entries of the reduced stiffness matrices share the same regularity as the coefficients. For the PDEs of interest, the coefficients have analytic regularity, so that the sparse approximations to $\widetilde{\mathbb{A}}_s^{\rm 00},\widetilde{\mathbb{A}}_s^{\rm 0 b},\widetilde{\mathbb{A}}_s^{\rm b0},\widetilde{\mathbb{A}}_s^{\rm bb}$ are expected to have spectral accuracy. In addition, perturbation theory is needed to analyze how such matrix approximation error propagates to PDE solutions.

\section{Numerical examples}\label{example}
To test the performance of our method, we carry out numerical experiments based on the two stochastic PDEs given in Example \ref{ex1} and \ref{ex2}, where the physical domain $D$ is set to a two-dimensional box $D := [0,1] \times [0,1]$. Our algorithms are implemented in \texttt{Matlab 2016a} and simulated on a workstation with \texttt{Intel(R) Xeon(R) CPU E5-2699 v4}. For each example, we will test two random fields. One is the truncated colored noise $\eta_N(x,\omega)$ defined in \eqref{e3} with the Gaussian covariance function 
\begin{equation}\label{cor_len}
\kappa(x,x') := \exp\left(-\frac{\|x-x'\|^2_2}{L^2}\right),
\end{equation}
where $L$ is the correlation length, and $\bm y := (y_1, \ldots, y_N)^{\top}$ are assumed to follow $N$-dimensional standard Gaussian distribution $\mathcal{N}(0, \mathbb{I})$. The other one is the discrete white noise $\eta^{\rm WN}_N(x, \omega)$ defined in \eqref{e4}, where the random variables  $y_1, \ldots, y_N$ are assumed to follow the $N$-dimensional Gaussian distribution $\mathcal{N}(0, \sigma^2 \mathbb{I})$ with $\sigma$ being the standard deviation. 

\subsection{Tests on the local KL expansion}
\begin{figure}[h!]
\center
\includegraphics[scale = 0.35]{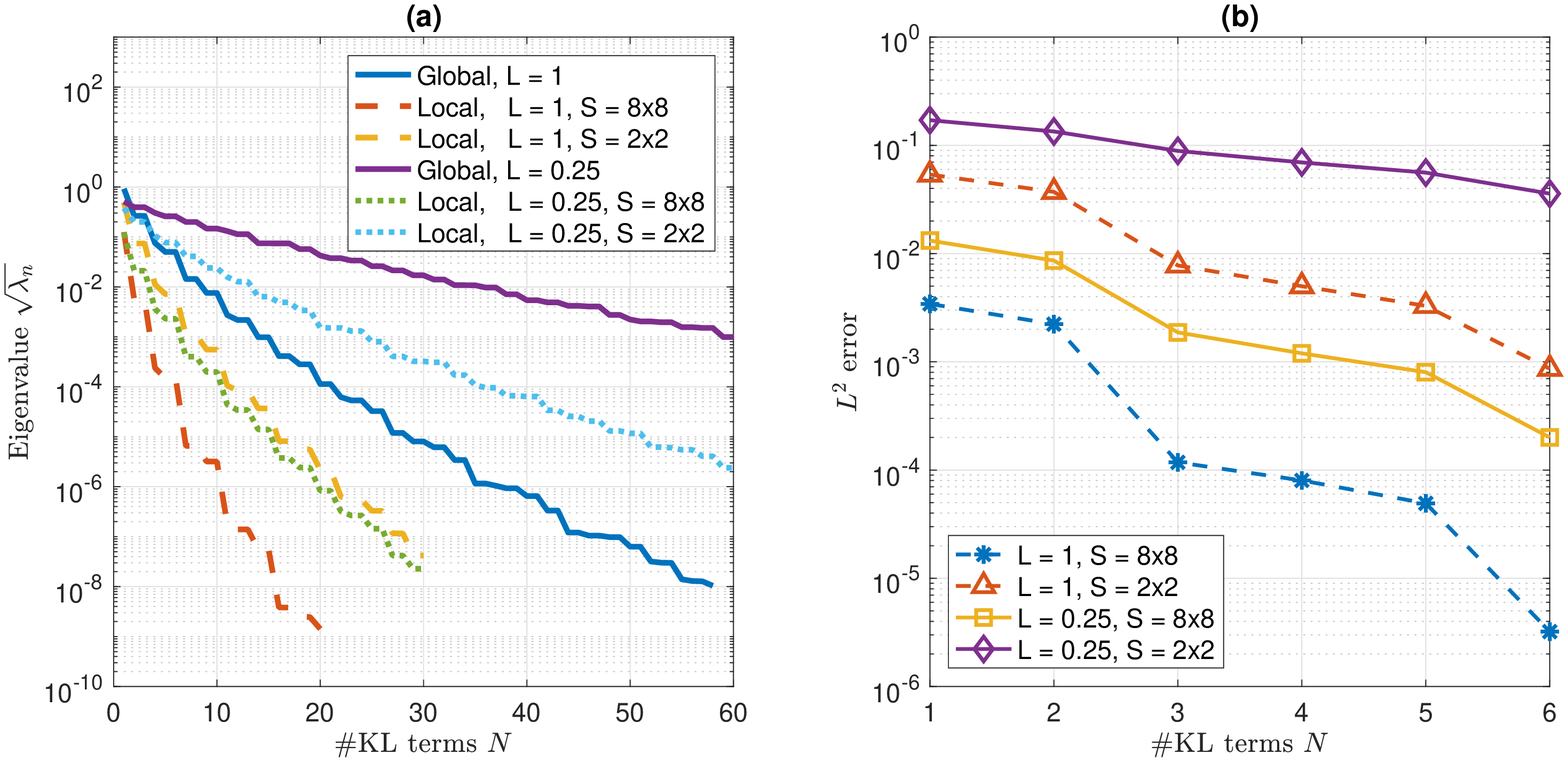}
\vspace{-0.3cm}
\caption{(a) The decay of eigenvalues $\sqrt{\lambda_n}$ of the global KL expansion $\eta_N$ in \eqref{e3} and the local KL expansion $\eta^{\rm loc}_{s,N_s}$ in \eqref{e11} for correlation length $L = 1$ and $L = 0.25$; (b) The $L^2$-error between the global and local truncated KL expansions, where the global KL expansion is truncated at the 200-th term.}\label{fig1}
\end{figure}

In the colored noise case, we need to test the error between the global and local KL expansions caused by the truncations and the least-squares projection, discussed in Section \ref{sec:locKL}. To do this, we define a $1024 \times 1024$ cartesian mesh in $D$, and discretize both the global and local random fields on the mesh $\mathcal{T}_h$ as piecewise constant functions. For simplicity, we assume the sub-domains are of the same size and shape, such that all the local KL expansions feature the same eigenvalue decay.  
In Figure \ref{fig1}(a), we show the decay of the eigenvalues $\sqrt{\lambda_n}$. As expected, the smaller the sub-domain, the faster the eigenvalues decay, which illustrates the motivation of using domain decomposition. In Figure \ref{fig1}(b), we compute the error $\|\eta_N - \eta_{\bm N}\|_{L^2(D)}$ for the four cases considered in Figure \ref{fig1}(a), where the global KL expansions are truncated at the 200-th term, and the local expansion are truncated at $N_s = 1,2,3,4,5,6$ for $s = 1, \ldots, S$. It can be seen that the error $\|\eta_N - \eta_{\bm N}\|_{L^2(D)}$ is dominated by the largest neglected eigenvalue of the local KL expansion. 

\subsection{The diffusion equation with random diffusivity}\label{ex:diff}
We consider the two-dimensional elliptic PDE given in Example \ref{ex1}, where $D = [0,1]^2$ and $f(x) = 100$. Piecewise linear finite element basis is used to discretize the PDE in $D$. 

\subsubsection{The colored noise case} \label{ex1:color}
\begin{figure}[h!]
\center
\includegraphics[scale = 0.25]{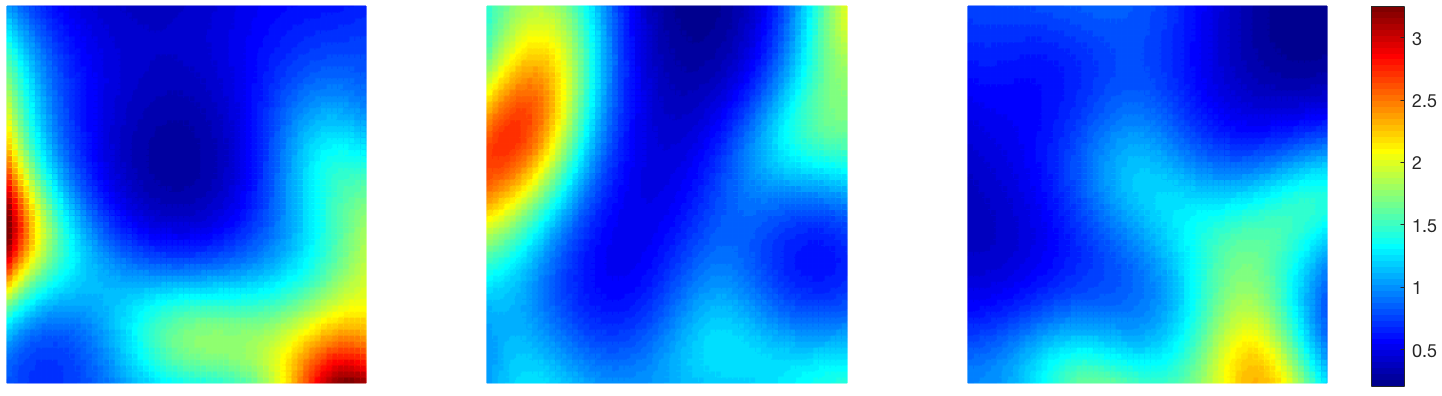}
\vspace{-0.3cm}
\caption{Illustration of three snapshots of the random field in \eqref{aaa} with $L = 0.25$.}\label{fig20}
\end{figure}

The diffusivity $a(x,\omega)$ in \eqref{eq:ellip} is defined by
\begin{equation}\label{aaa}
a(x,\omega) := \exp\left( \frac{1}{5}\eta(x,\omega)\right),
\end{equation}
with $\eta(x, \omega)$ defined in \eqref{e10} based on the covariance function in \eqref{cor_len} with the correlation length being $L = 1$ and $L = 0.25$. Figure \ref{fig20} shows three snapshots of the random field in \eqref{aaa} with $L = 0.25$. To compute the total approximation error, we define the reference solution to be the numerical solution obtained by solving the PDE in \eqref{eq:ellip} on a triangle mesh with $h = 1/2^{12}$, i.e., a total of $4096 \times 4096$ grid points, using the truncated global KL expansion $\eta_N$ with $N = 200$.
Figure \ref{fig1}(a) shows that $N= 200$ is sufficient to neglect the global KL truncation error. 
The random variables $y_1, \ldots, y_N$ follow the $N$-dimensional standard Gaussian distribution. For the local KL expansion, 
the bounded domains $\Gamma_s$, introduced in Section \ref{sec:data}, are set to $\Gamma_s = [-5,5]^{N_s}$ for $s = 1, \ldots, S$, such that the probability of having a sample $\bm y_s$ fall outside $\Gamma_s$ is about $5 \times 10^{-7}$. In the domain $\Gamma_s$, we use sparse Legendre polynomials to approximate the local stiffness matrices. 

\begin{figure}[h!]
\center
\includegraphics[scale = 0.35]{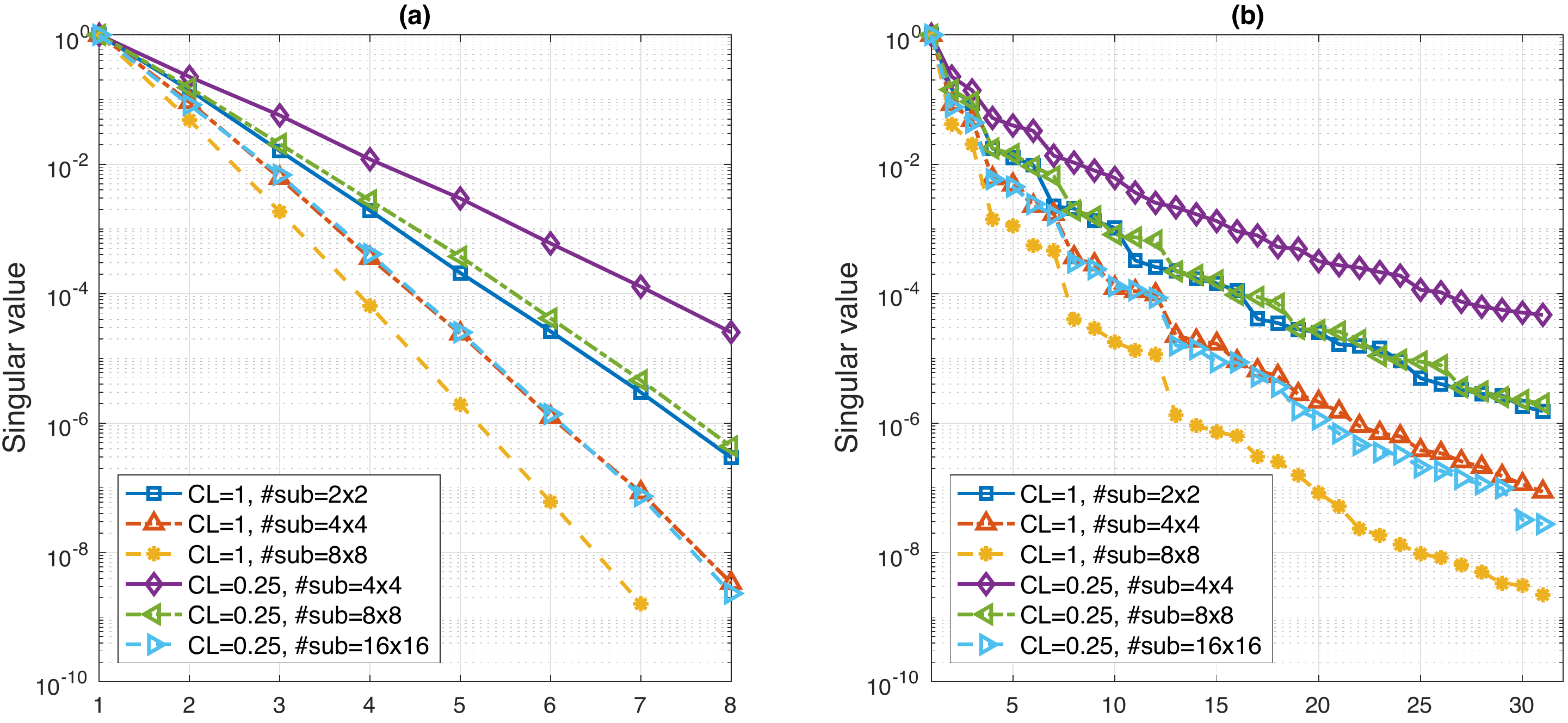}
\vspace{-0.3cm}
\caption{(a) The decay of singular values of the solution on one interface; (b) The decay of singular values of the solution in the interior of one sub-domain $D_s^H$.}\label{fig2}
\end{figure}

To illustrate the effectiveness of SVD, we draw 1000 random samples of $\bm y \in \mathbb{R}^N$, generate snapshots by solving the expensive finite element problems, perform SVDs, and plot in Figure \ref{fig2} the singular value\footnote{The plotted singular values are normalized by the largest singular value in each case.} decays along an interface and in the interior of a sub-domain. We observe that the more sub-domains, the faster the singular values decay. As such, we can keep a small number of singular vectors along the interfaces to reduce the size of the global stiffness matrix $\mathbb{B}$ in \eqref{globalsys}, as well as keep a small number of interior singular vectors in each sub-domain to reduce the size of the local stiffness matrices $\mathbb{A}_s^{00}$ in \eqref{schur}. Figure \ref{fig2} demonstrates that obtaining fast singular value decay is another advantage of using domain decomposition. 

Now we show the accuracy of our DDMR approach by examining the error decay of the reduced model with respect to 4 quantities, i.e., 
\begin{itemize}
\item $N_s$: the dimension  of the truncated local KL expansion in \eqref{e11};
\item $M_{s,j}$: the dimension  of $\widehat{\mathbb{V}}_{s,j}^{\rm b}$ in \eqref{sing} along the $j$-th interface of $D_s^H$;
\item $M_s$: the dimension  of $\widehat{\mathbb{V}}_s^{0}$ in \eqref{sing1} in sub-domain $D_s^H$;
\item $M^{\rm LS}_{s,p}$: the cardinality of the polynomial space $\mathcal{P}_{M_{s,p}^{\rm LS}}(\Gamma_s)$ of \eqref{DLS}.
\end{itemize}
Each of the above four quantities could be different for each sub-domain or interface. In this work, we will restrict us to use the same number over all sub-domains for each quantity. This strategy is not optimal, but sufficient to demonstrate the performance of our method. 
We run simulations in 4 different scenarios, i.e.,
\begin{itemize}
\item[(i)] $L = 1$ with $S = 8 \times 8 $ sub-domains;
\item[(ii)] $L = 1$ with $S = 16 \times 16$ sub-domains;
\item[(iii)] $L = 0.25$ with $S = 8 \times 8 $ sub-domains;
\item[(iv)] $L = 0.25$ with $S = 16 \times 16$ sub-domains.
\end{itemize}
\begin{figure}[h!]
\center
\includegraphics[scale = 0.4]{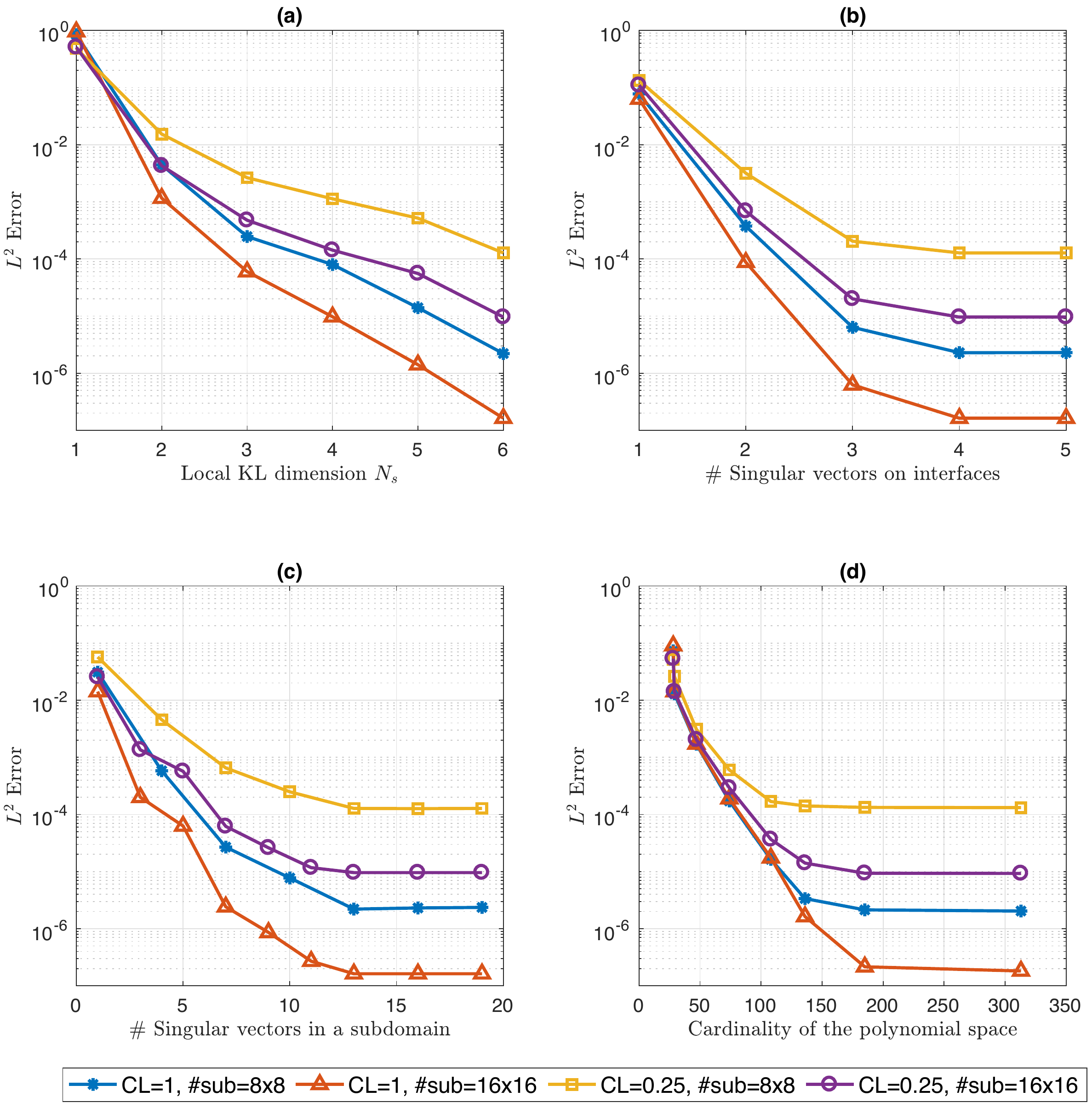}
\vspace{-0.2cm}
\caption{Illustration of the error decays w.r.t.~$N_s$, $M_{s,j}$, $M_s$ and $M_{s,p}^{\rm LS}$. (a) Error decay w.r.t.~$N_s$ while setting $M_{s,j} = 6$, $M_s = 19$ and $M_{s,p}^{\rm LS} = 332$; (b) Error decay w.r.t.~$M_{s,j}$ while setting $ N_s = 6$, $M_s = 19$ and $M_{s,p}^{\rm LS} = 332$; (c) Error decay w.r.t.~$M_{s}$ while setting $ N_s = 6$, $M_{s,j}= 6$ and $M_{s,p}^{\rm LS} = 332$; (d) Error decay w.r.t.~$M_{s,p}^{\rm LS}$ while setting $ N_s = 6$, $M_{s,j}= 6$ and $M_s = 19$.}\label{fig3}
\end{figure}
The results are shown in Figure \ref{fig3}. 
In each scenario, we first compute the errors by setting $N_s = 6$, $M_{s,j} = 6$, $M_{s} = 19$ and $p = 9$ (leading to $M_{s,p}^{\rm LS} = 332$). These errors are shown as the last markers of the error decay curves plotted in Figure \ref{fig3}.
Then, in each sub-figure of Figure \ref{fig3}, we decrease one of the four quantities while remaining the other three unchanged.
In each scenario, we generate 1000 random samples of $\bm y$ and execute 1000 time-consuming finite element solvers (using the $4096 \times 4096$ mesh) to generate data of the solution defined in \eqref{data}. On the other hand, we also generate another 1000 realizations of $\bm y_s$ for $s = 1, \ldots, S$ to generate the data of the local stiffness matrices defined in \eqref{appAF}. The error is computed in the relative $L^2$ norm using another 1000 random samples of $\bm y \sim\mathcal{N}(0, \mathbb{I})$. The error decays as what we expected.
 We would like to point out that the smallest error comes from the case of having a bigger correlation and more sub-domains. Thus, for a small correlation length, we can introduce more sub-domains to reduce the total error.

%

Next, we discuss the efficiency of our approach to achieve a prescribed accuracy. We use the case with $L = 0.25$ and three triangle meshes of sizes $1024 \times 1024$, $2048 \times 2048$ and $4096 \times 4096$, respectively. For each of the three meshes, we construct our reduced model by setting $N_s = 6$, $M_{s, j} = 6$, $M_{s} = 19$ and $p= 9$ (i.e., $M_{s,p}^{\rm LS} = 332$) in our algorithm\footnote{The model reduction error is balanced with the FE error on the finest mesh of size $4096\times 4096$.}. The complexity of our approach is divided into the offline cost, i.e., the cost of Algorithm 1, and the online cost, i.e., the cost of Algorithm 2. Both the online and the offline costs are measured in a relative way by 
\begin{equation}\label{cost}
\rm Cost = \frac{\text{The offline(online) CPU time}}{\text{The CPU time of one expensive FE simulation}},
\end{equation}
which, in other words, is the number of expensive FE simulations. To be more harsh to our method, the CPU time of one expensive FE simulation is only measured by the time of solving the final linear system using the ``$\backslash$'' solver in \texttt{Matlab 2016a}, regardless of assembly cost and other operations that may not be optimally implemented.
\setlength{\extrarowheight}{3pt}
\begin{table}[h!]
\footnotesize
\caption{Computational cost for solving Example \ref{ex1} with the diffusivity being the colored noise ($L = 0.25$). The domain $D$ is decomposed into $16\times 16$ sub-domains. 
The unit FE time only takes into account the CPU time of solving the final linear system; the offline and online costs are measured by the number of expensive FE simulations.}\label{ex1_t1}
\begin{center}
\begin{tabular}{c|c|c|c|c}
\hline
\multirow{2}{*}{FE cost}& \# FE  nodal values & $2^{20}$ & $2^{22}$ & $2^{24}$\\
\cline{2-5}
& Unit FE time & 4 sec & 25 sec &	381 sec \\
\hline\hline
\multirow{6}{*}{\shortstack{Offline cost \\ \big($\frac{\text{Wall time}}{\text{Unit FE time}}$\big)}}& KL expansion &    2.19  & 2.16 & 	0.53\\
& FE solves for $\mathbf{U}_s(w_k), k = 1, \ldots, 1000$ & 1000 & 1000 & 1000\\
& SVD on the interfaces for $\widehat{\mathbb{V}}_{s,j}^{\rm b}$ &  22.70 & 4.71 &	0.33\\
& SVD in the sub-domains for $\widehat{\mathbb{V}}_{s}^{\rm 0}$ & 34.13 & 29.00 & 12.38\\
& Assembling $\widetilde{\mathbb{A}}_s(\omega_k), k = 1, \ldots, 1000$ & 256.25 & 301.33 & 389.78\\
& Computing $\widehat{\mathbb{K}}_s^{00}, \widehat{\mathbb{K}}_s^{\rm 0b}, \widehat{\mathbb{K}}_s^{\rm b0}, \widehat{\mathbb{K}}_s^{\rm bb}$ & 6.89	& 1.45 & 0.08\\ \cline{2-5}
& Total & 1315.27 & 1338.65 & 1403.1\\
\hline\hline
\multirow{3}{*}{\shortstack{Online cost\\ \big($\frac{\text{Wall time}}{\text{Unit FE time}}$\big)}}& Solving $ \widehat{\mathbf{C}}^{\rm b, LS}$ &  0.09  & 0.02 &  0.001 \\
& Solving  $\widehat{\mathbf{C}}^{\rm 0, LS}_s$ for $s = 1, \ldots, S$ & 0.33 & 0.05 & 0.003\\
\cline{2-5}
& Total & 0.42 & 0.07 & 0.004\\
\hline
\end{tabular}\vspace{-0.0cm}
\end{center}
\end{table}
The results are shown in Table \ref{ex1_t1} for the case of having $S = 16 \times 16$ sub-domains. In the offline procedure, the data generation is still the dominant part, and our dimension reduction strategy successfully helps achieve $\mathcal{O}(10^{-6})$ error with only 1000 expensive FE simulations. For the SVD algorithm, since we only need a few of the largest singular values and singular vectors, we do not need to run full SVD. Instead, we use the Lanczos bi-diagonalization methods \cite{Baglama:2005ji} to reduce the cost of running SVDs. On the other hand, the wall time of the online cost is $0.42\times 4\sec = 1.68\sec$, $0.07\times 25\sec = 1.75\sec$, $0.004\times381\sec = 1.52\sec$ for the cases of having $1024 \times 1024$, $2048 \times 2048$ and $4096 \times 4096$ meshes, respectively. This verifies that the online cost is independent of the triangle mesh size $h$, so that the finer the original mesh, the more savings our method can provide.

%
%

\subsubsection{The discrete white noise case}\label{sec:ex21}
\begin{figure}[h!]
\center
\includegraphics[scale = 0.53]{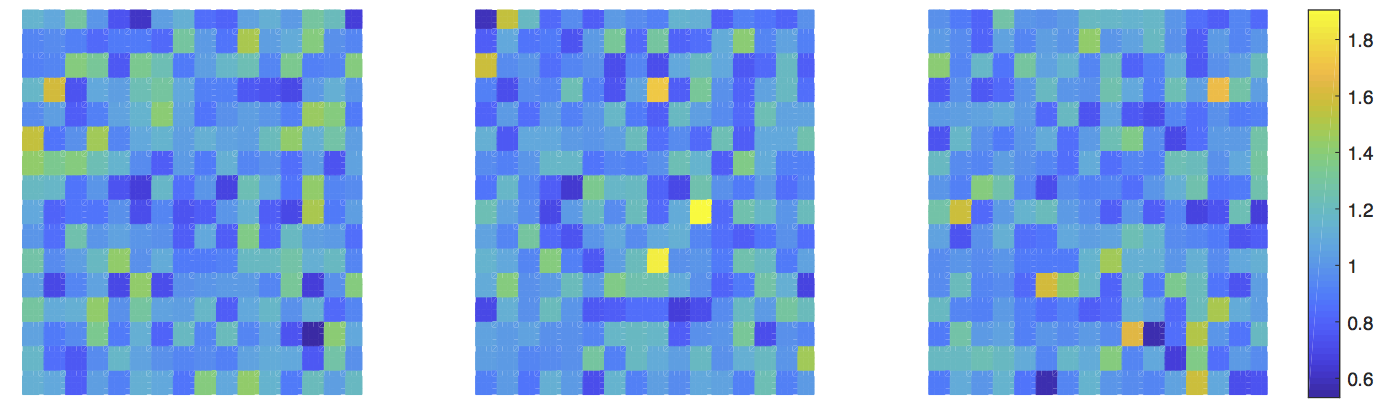}
\vspace{-0.2cm}
\caption{Illustration of three snapshots of the discrete white noise with $N = 16 \times 16$ and $\sigma = 1$.}\label{fig44}
\end{figure}

Now we test the discrete white noise case by replacing the random field $\eta$ in \eqref{aaa} with the one in \eqref{e4}, which is assumed to be a uniformly partitioned piecewise constant random field.
The random parameters $\bm y$ follow the multi-dimensional Gaussian distributions, denoted by $\mathcal{N}(0, \sigma^2 \mathbb{I})$, where $\sigma$ is the standard deviation and $\mathbb{I}$ is the $N$-dimensional identity matrix. Figure \ref{fig44} shows three snapshots of the discrete white noise with $N = 16 \times 16$ and $\sigma = 1$.

We run simulations in 4 different cases, i.e.,
\begin{itemize}
\item[(i)] $\sigma = 0.1$ with $N = S= 8 \times 8 $ sub-domains;
\item[(ii)] $\sigma = 0.1$ with $N =S= 16 \times 16$ sub-domains;
\item[(iii)] $\sigma = 1$ with $N =S= 8 \times 8 $ sub-domains;
\item[(iv)] $\sigma = 1$ with $N =S= 16 \times 16$ sub-domains.
\end{itemize}
For each case, we assume $D_n^{\rm WN} = D_s^{H}$, i.e., align the domain decomposition with the partition of the random field. As such, the approximation of the local stiffness matrices becomes a one-dimensional approximation problem. We still use 1000 realizations to conduct SVD, but we only use 100 realizations for the DLS approximation of the local stiffness matrices. The local parameter domain $\Gamma_s$ in \eqref{sec:data} is set to $[-5\sigma, 5\sigma]$. The reference solution is obtained by solving the PDE on a very fine mesh with $h = 1/2^{12}$, i.e., $4096 \times 4096$ unknowns. 
All the other settings are the same as in Section \ref{ex1:color}. 
\begin{figure}[h!]
\center
\includegraphics[scale = 0.4]{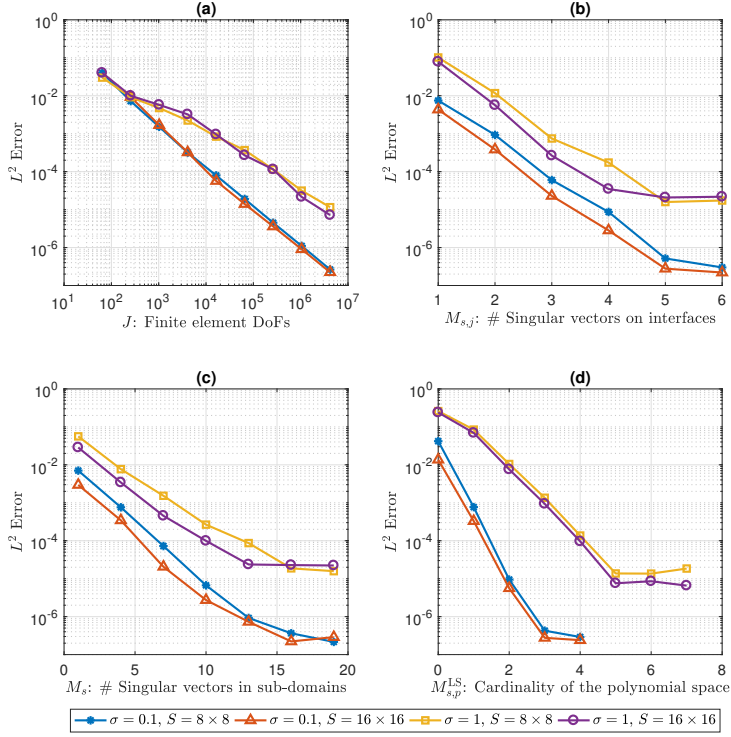}
\vspace{-0.2cm}
\caption{Illustration of the error decays w.r.t.~$J$, $M_{s,j}$, $M_s$ and $M_{s,p}^{\rm LS}$. (a) Error decay w.r.t.~$J$ while setting $M_{s,j} = 6$, $M_s = 19$ and $M_{s,p}^{\rm LS} = 10$; (b) Error decay w.r.t.~$M_{s,j}$ while setting $ J = 2^{24}$, $M_s = 19$ and $M_{s,p}^{\rm LS} = 10$; (c) Error decay w.r.t.~$M_{s}$ while setting $ J = 2^{24}$, $M_{s,j}= 6$ and $M_{s,p}^{\rm LS} = 10$; (d) Error decay w.r.t.~$M_{s,p}^{\rm LS}$ while setting $J = 2^{24}$, $M_{s,j}= 6$ and $M_s = 19$.}\label{fig4}
\end{figure}

The results are shown in Figure \ref{fig4}, where we use the same strategy as in Figure \ref{fig3} to generate the error decay curves. Since there is no KL expansion setting, we plotted the error with respect to $J$ in \eqref{eq:FEM}. The error is computed in relative $L^2$ norm using another 1000 random samples of $\bm y \sim\mathcal{N}(0, \sigma^2 \mathbb{I})$. As expected, the error of the case (ii) is the smallest due to the smaller variance. The flat toes in Figure \ref{fig4}(b)-(d) are due to the dominancy of the error caused by the triangle mesh size. As shown in Figure \ref{fig4}(d), an important advantage of our method is that it can reduce local stiffness matrix approximation to a set of one-dimensional problems that completely overcomes the curse of dimensionality. 

\setlength{\extrarowheight}{3pt}
\begin{table}[h!]
\footnotesize
\caption{Computational cost for solving Example \ref{ex1} with diffusivity being the 256-dimensional discrete white noise. The standard deviation of each random variable is $\sigma = 0.1$. The unit FE time only takes into account the CPU time of solving the final linear system; the offline and online costs are measured by the number of expensive FE simulartions.}\label{ex1_t3}
\begin{center}
\begin{tabular}{c|c|c|c|c}
\hline
\multirow{2}{*}{
FE cost}& \# FE  nodal values & $2^{20}$ & $2^{22}$ & $2^{24}$\\
\cline{2-5}
& Unit FE time & 4 sec & 27 sec &	394 sec \\
\hline\hline
\multirow{5}{*}{\shortstack{Offline cost \\ \big($\frac{\text{Wall time}}{\text{Unit FE time}}$\big)}}
&  FE solves for $\mathbf{U}_s(w_k), k = 1, \ldots, 1000$ & 1000 & 1000	& 1000\\
& SVD on the interfaces for $\widehat{\mathbb{V}}_{s,j}^{\rm b}$& 23.82 & 4.17& 0.32\\
& SVD in the sub-domains for $\widehat{\mathbb{V}}_s^0$ & 35.42 & 29.51& 12.18\\
& Assembling $\widetilde{\mathbb{A}}_s(\omega_k), k = 1, \ldots, 100$ & 28.90 & 34.83 & 49.67\\
& Computing $\widehat{\mathbb{K}}_s^{00}, \widehat{\mathbb{K}}_s^{\rm 0b}, \widehat{\mathbb{K}}_s^{\rm b0}, \widehat{\mathbb{K}}_s^{\rm bb}$ & 0.005 &	 0.007 &	0.001\\ \cline{2-5}
& Total & 1088.15 &	1068.52 & 1062.17\\
\hline\hline
\multirow{3}{*}{\shortstack{Online cost \\ \big($\frac{\text{Wall time}}{\text{Unit FE time}}$\big)}}& Solving $ \widehat{\mathbf{C}}^{\rm b, LS}$ &  0.111	& 0.018 &	0.001 \\
& Solving  $\widehat{\mathbf{C}}^{\rm 0, LS}_s$ for $s = 1, \ldots, S$& 0.423 & 0.071 &	0.002\\
\cline{2-5}
& Total & 0.534 &	0.089&	0.003\\
\hline
\end{tabular}\vspace{-0.0cm}
\end{center}
\end{table}
Next, we discuss the efficiency of our approach using the 256-dimensional random field, i.e., $N = S = 16 \times 16$ in \eqref{e4}, with $\sigma = 0.1$. Three triangle meshes of sizes $1024 \times 1024$, $2048 \times 2048$ and $4096 \times 4096$ are used. For each mesh, we construct a reduced model by setting $M_{s, j} = 6$, $M_{s} = 19$ and $p= 9$ (i.e., $M_{s,p}^{\rm LS} = 10$) in our algorithm. The complexity of our approach is divided into the offline cost, i.e., the cost of Algorithm 1, and the online cost, i.e., the cost of Algorithm 2. Both the online and the offline costs are measured in a relative way by the formula given in \eqref{cost}. The results are shown in Table \ref{ex1_t3}. As expected, our managed to make the online cost independent of the original triangle mesh size. In the offline procedure, the cost of assembling $\widetilde{A}_s$ is much smaller than the case of have colored noise, because of the one-dimensional parametric dependence of the local stiffness matrices.

\subsection{The convection-dominated transport with random velocity}\label{CD}
We consider the two-dimensional PDE given in Example \ref{ex2}, where $f(x) = 0$, $D = [0,1]^2$ and the set $\mathscr{D}$ is a subset of $\partial D$ defined by $\{x_1 = 0, x_2 \in [0,0.5]\} \cup \{x_1 \in [0,1], x_2 =0\}$. The SUPG scheme is used to discretize the PDE in the physical domain $D$. 

\subsubsection{The colored noise case} \label{ex2:color}
\begin{figure}[h!]
\center
\includegraphics[scale = 0.4]{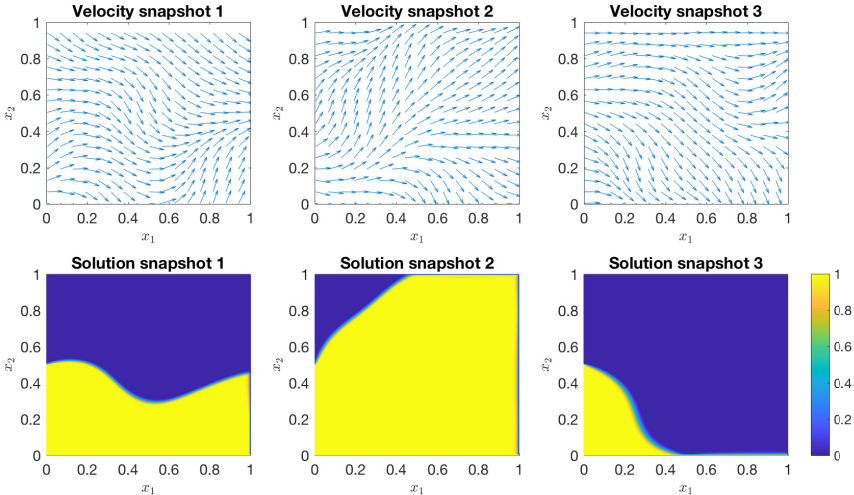}
\vspace{-0.5cm}
\caption{Three snapshots of the colored noise velocity field defined in \eqref{bbb} and the corresponding solution field of the PDE in \eqref{eq:conv}}\label{fig7}
\end{figure}

The random velocity field $\bm b(x,\omega)$ is defined by 
\begin{equation}\label{bbb}
\bm b(x,\omega) := 
\begin{pmatrix}
\cos\left(\dfrac{1}{5} \eta(x, \omega)\right)  \vspace{0.2cm}\\
\sin\left(\dfrac{1}{5} \eta(x, \omega) \right)
\end{pmatrix},
\end{equation}
where $\eta(x,\omega)$ is defined in \eqref{e10} based on the covariance function given in \eqref{cor_len} with correlation length $L = 0.25$. The reference solution is obtained by solving the PDE in \eqref{eq:conv} on a mesh with $h = 1/2^{12}$ using the truncated global KL expansion $\eta_N$ with $N = 200$. The random variables $y_1, \ldots, y_N$ follow the $N$-dimensional standard Gaussian distribution $\mathcal{N}(0, \mathbb{I})$. Three snapshots of the velocity field and the corresponding solution field are given in Figure \ref{fig7}. 

We define a total of $S = 16 \times 16$ sub-domains.
 The bounded domain $\Gamma_s$ introduced in Section \ref{sec:data} are set to $\Gamma_s = [-5,5]^{N_s}$ for $s = 1, \ldots, S$, as in Example \ref{ex1}. 
Similar to Figure \ref{fig3}, we show the accuracy of our approach by examining the error decays of the reduced model with respect to the 4 quantities, i.e., $N_s$, $M_{s,j}$, $M_s$ and $M^{\rm LS}_{s,p}$. 
We run simulations with 3 diffusion coefficients, i.e., $\varepsilon  = 10^0, 10^{-2}, 10^{-4}$.
\begin{figure}[h!]
\center
\includegraphics[scale = 0.4]{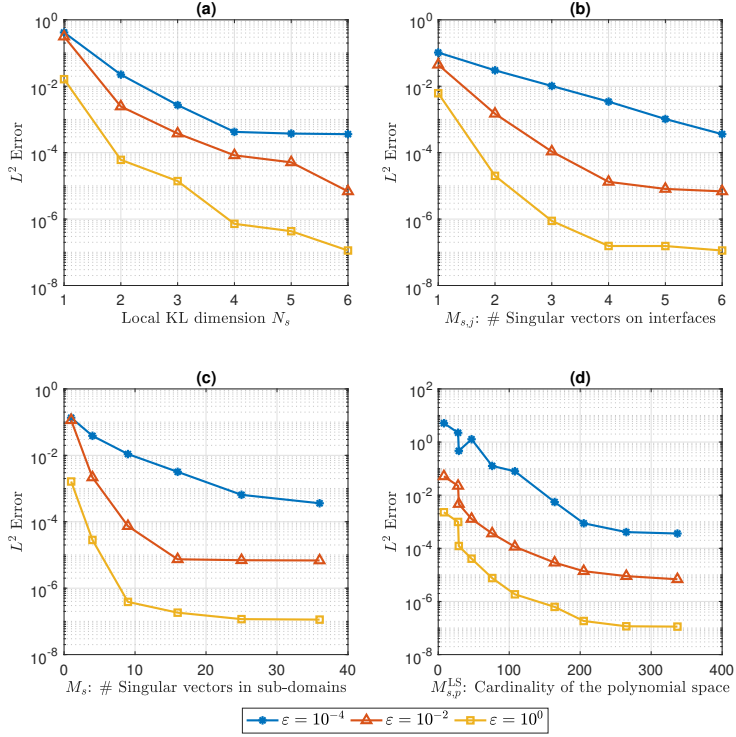}
\vspace{-0.2cm}
\caption{Illustration of the error decays w.r.t.~$N_s$, $M_{s,j}$, $M_s$ and $M_{s,p}^{\rm LS}$. (a) Error decay w.r.t.~$N_s$ while setting $M_{s,j} = 6$, $M_s = 36$ and $M_{s,p}^{\rm LS} = 332$; (b) Error decay w.r.t.~$M_{s,j}$ while setting $ N_s = 6$, $M_s = 36$ and $M_{s,p}^{\rm LS} = 332$; (c) Error decay w.r.t.~$M_{s}$ while setting $ N_s = 6$, $M_{s,j}= 6$ and $M_{s,p}^{\rm LS} = 332$; (d) Error decay w.r.t.~$M_{s,p}^{\rm LS}$ while setting $ N_s = 6$, $M_{s,j}= 6$ and $M_s = 36$.}\label{fig5}
\end{figure}

The results are shown in Figure \ref{fig7}. 
We first compute the errors by setting $N_s = 6$, $M_{s,j} = 6$, $M_{s} = 36$ and $p = 9$ (leading to $M_{s,p}^{\rm LS} = 332$). These errors are shown as the last markers of the error decay curves plotted in Figure \ref{fig7}. Then we use the same strategy as in Figure \ref{fig3} to generate
Figure \ref{fig7}.
%
We observe that the errors become bigger as we decrease the value of $\varepsilon$, i.e., increasing the sharpness of the transition area. Moreover, when $\varepsilon = 10^0$, the solution has a smooth transition, in which case the error decays fast with respect to $M_{s,j}$ and $M_s$. In comparison, when $\varepsilon = 10^{-4}$, the error from SVDs are dominant, as both the errors caused by local KL expansion and the sparse approximation hit flat toes. In addition, we also observe that the decay rate of the error with respect to $M_{s,p}^{\rm LS}$ remains the same as we decrease the value of $\varepsilon$, because of the fact that the smoothness of the local stiffness matrices is not changed with $\varepsilon$.

\subsubsection{The discrete white noise case}
\begin{figure}[h!]
\center
\includegraphics[scale = 0.4]{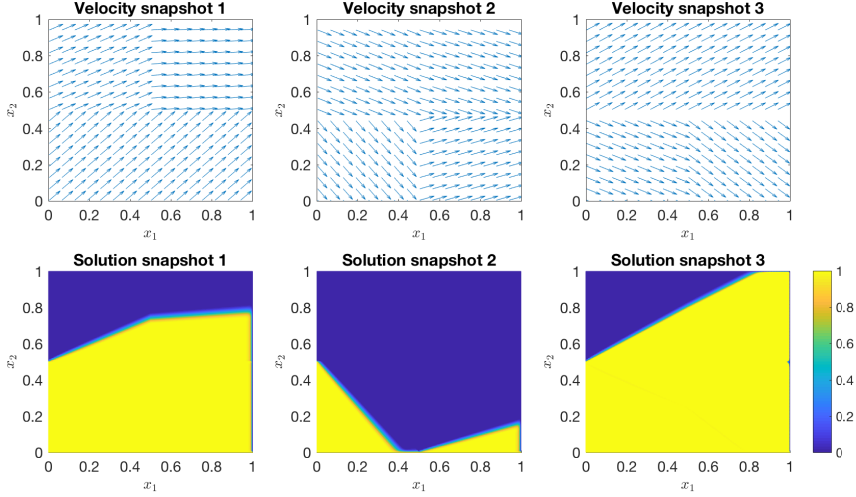}
\vspace{-0.5cm}
\caption{Three snapshots of the discrete white noise velocity field defined in \eqref{bbb} and the corresponding solution field of the PDE in \eqref{eq:conv}.}\label{fig8}
\end{figure}

Now we test the discrete white noise case by replacing the random field $\eta$ in \eqref{bbb} with the one in \eqref{e4}, which is a uniformly partitioned piecewise constant random field.
The random parameters in $\bm y$ follow $\mathcal{N}(0, \sigma^2 \mathbb{I})$ as in Section \ref{sec:ex21}.
In the case of
$N = 2 \times 2$, $\sigma = 0.5$, three snapshots of the velocity field and the corresponding solution field are given in Figure \ref{fig8} for illustration.

To test the accuracy, we run simulations in 4 different cases, i.e.,
\begin{itemize}
\item[(i)] $\sigma = 0.5$ with $\varepsilon = 10^{-4}$;
\item[(ii)] $\sigma = 0.5$ with $\varepsilon = 10^{-2}$;
\item[(iii)] $\sigma = 0.1$ with $\varepsilon = 10^{-4}$;
\item[(iv)] $\sigma = 0.1$ with $\varepsilon = 10^{-2}$.
\end{itemize}
For each case, we assume $D_n^{\rm WN} = D_s^{H}$ as in Example \ref{ex1}. Again, the approximation of the reduced local stiffness matrices becomes a one-dimensional approximation problem. We still use 1000 realizations to conduct SVD, but we only use 100 realizations for the DLS approximation of the local stiffness matrices. The local parameter domain $\Gamma_s$ in \eqref{sec:data} is set to $[-5,5]$. 
All the other settings are the same as in Section \ref{sec:ex21}. 
\begin{figure}[h!]
\center
\includegraphics[scale = 0.3]{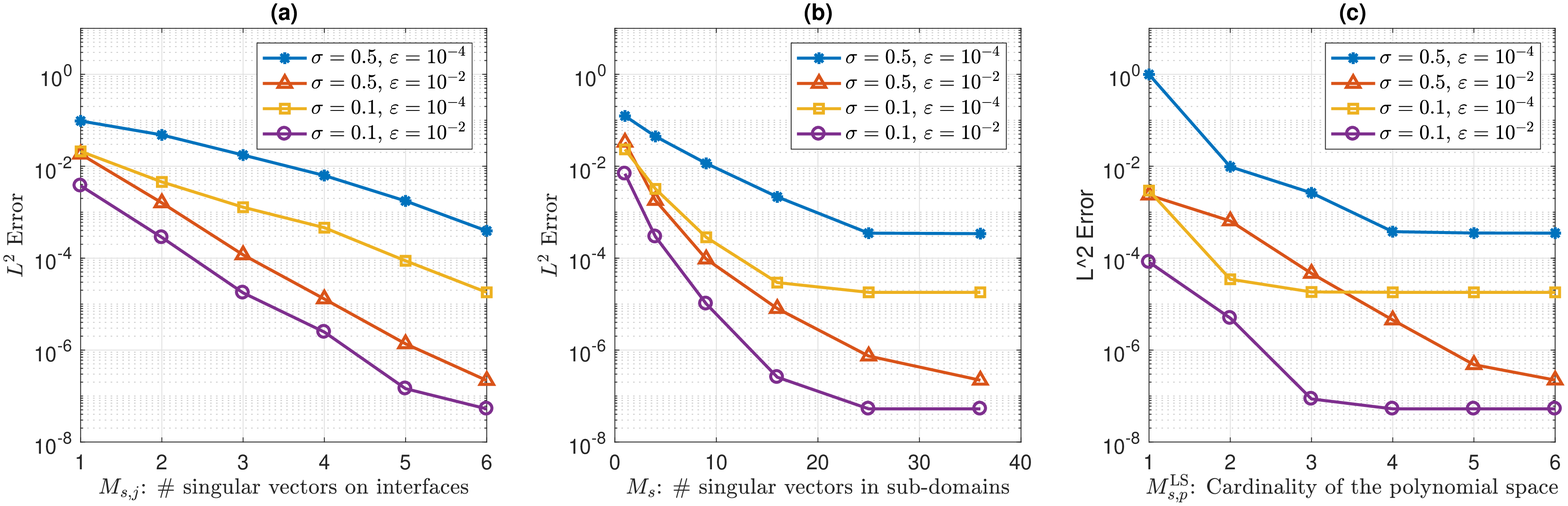}
\vspace{-0.2cm}
\caption{Illustration of the error decays w.r.t.~$M_{s,j}$, $M_s$ and $M_{s,p}^{\rm LS}$. (a) Error decay w.r.t.~$M_{s,j}$ while setting $M_s = 36$ and $M_{s,p}^{\rm LS} = 7$; (b) Error decay w.r.t.~$M_{s}$ while setting $M_{s,j}= 6$ and $M_{s,p}^{\rm LS} = 7$; (c) Error decay w.r.t.~$M_{s,p}^{\rm LS}$ while setting $M_{s,j}= 6$ and $M_s = 36$.}\label{fig6}
\end{figure}
The results are shown in Figure \ref{fig6}. The error is computed in relative $L^2$ norm using another 1000 random samples of $\bm y \sim \mathcal{N}(0, \sigma^2 \mathbb{I})$. As expected, the error of the case (iv) is the smallest due to the smaller $\sigma$ and the bigger $\varepsilon$. As $\varepsilon$ decreases, the SVD errors become more and more dominant, due to the slow decay of singular values around the sharp transition layer. 

\begin{remark}
For computational cost of our method, we can conduct similar discussions and draw the same conclusions as in Example \ref{ex1} by generating tables analogous to Table \ref{ex1_t1} and \ref{ex1_t3}. Thus, we omit the cost analysis in this example and only refer to the discussions on complexity in Section \ref{ex:diff}. 
\end{remark}

\section{Concluding remarks}\label{conc}
We developed a new model reduction method for stochastic convection-diffusion equations by integrating domain decomposition, local reduced basis method, and sparse approximation of operators. Our method can overcome the curse of high-dimensionality, achieves online-offline decomposition, as well as handle the convection-dominated problem with irregular behavior. Even though the our strategy shows very promising performance, it could be further improved in several aspects. The first direction would be incorporating greedy algorithms \cite{DeVore:2013eq,Hesthaven:2014er} to generate the local reduced bases, which requires a-posteriori error estimates and a systematic way to coordinate the greedy search in different sub-domains. Moreover, more advanced Galerkin formulations, e.g., the weak Galerkin methods \cite{Mu:2013fb,Mu:2017fy}, could be incorporated to extend our method for other type of PDEs. 
Another direction is to extend our method to a non-intrusive algorithm (e.g., \cite{Peherstorfer:2016je}), which will make it much easier to couple with large-scale simulation code. 

\section*{Acknowledgement}
We thank Professor Christoph Schwab for his valuable comments and suggestions on this effort and future research directions on this topic. 

\bibliographystyle{siam}


\end{document}